\documentclass[bj]{imsart}
\RequirePackage{amsthm,amsmath,amsfonts,amssymb}
\RequirePackage[numbers]{natbib}
\RequirePackage[colorlinks,citecolor=blue,urlcolor=blue]{hyperref}
\RequirePackage{graphicx}

\usepackage{eucal,xr}
\usepackage{amssymb,epsfig,latexsym,xcolor}
\usepackage{graphicx,amsmath,amssymb,latexsym,psfrag}
\usepackage{graphics,graphicx,comment,bm,empheq}
\usepackage{multirow,changepage}
\usepackage[subnum]{cases}
\startlocaldefs

\newcommand{\cG}{\mathcal{G}}
\newcommand{\cA}{\mathcal{A}}

\newcommand{\cU}{\mathcal{U}}

\newcommand{\N}{\mathbb{N}}

\newtheorem{theorem}{Theorem}[section]
\newtheorem{lemma}[theorem]{Lemma}

\newtheorem{proposition}[theorem]{Proposition}
\theoremstyle{remark}

\newtheorem{remark}[theorem]{Remark}

\newcommand{\oD}{\mathcal{D}}

\newcommand{\oE}{\mathcal{E}}

\newcommand{\osn}[1]{n_{#1}}

\newcommand{\oboldm}{\bold{u}}

\newcommand{\oldb}{\rho}

\newcommand{\sidebyside}[4]{
\begin{figure}
\begin{minipage}[t]{6 cm}
\includegraphics[width=0.80\textwidth]{#1.pdf}
\caption{#2\label{fig:#1}}
\end{minipage}
\hfill \hspace{0.4cm}
\begin{minipage}[t]{6 cm}
\includegraphics[width=0.80\textwidth]{#3.pdf}
\caption{#4\label{fig:#3}}
\end{minipage}
\hfill
\vspace{-0mm}
\end{figure}
}

\newcommand{\threplots}[6]{
\begin{figure}
\begin{minipage}[t]{4.4 cm}
\includegraphics[width=0.98\textwidth]{#1.pdf}
\caption{#2\label{fig:#1}}
\end{minipage}
\hfill \hspace{0.2cm}
\begin{minipage}[t]{4.4 cm}
\includegraphics[width=0.98\textwidth]{#3.pdf}
\caption{#4\label{fig:#3}}
\end{minipage}
\hfill \hspace{0.2cm}
\begin{minipage}[t]{4.4 cm}
\includegraphics[width=0.98\textwidth]{#5.pdf}
\caption{#6\label{fig:#5}}
\end{minipage}
\hfill
\vspace{-0mm}
\end{figure}
}
\endlocaldefs

\begin{document}
	
		\begin{frontmatter}
		\title{Bootstrap percolation on the stochastic block model}
		\runtitle{Bootstrap percolation on the SBM}
		
		\begin{aug}
			\author[A]{\fnms{Giovanni Luca} \snm{Torrisi}\ead[label=e1]{giovanniluca.torrisi@cnr.it}},
			\author[B]{\fnms{Michele} \snm{Garetto}\ead[label=e2]{michele.garetto@unito.it}}
			\and
			\author[C]{\fnms{Emilio} \snm{Leonardi}\ead[label=e3]{emilio.leonardi@polito.it}}
			\address[A]{CNR-IAC, Roma, Italy.
				\printead{e1}}
			
			\address[B]{Universit\`a di Torino, Torino, Italy.
				\printead{e2}}
		
			\address[C]{Politecnico di Torino, Torino, Italy.
			\printead{e3}}
		\end{aug}



\begin{abstract}
We analyze the bootstrap percolation process on the stochastic
block model (SBM), a natural extension of the Erd\H{o}s--R\'{e}nyi random graph
that incorporates the community structure observed in many real systems.
In the SBM, nodes are partitioned into two subsets, which represent different communities,
and pairs of nodes are independently connected with a probability that depends on the communities they belong to.
Under mild assumptions on the system parameters, we prove the existence of a sharp phase transition
for the final number of active nodes and characterize the sub-critical and the super-critical regimes
in terms of the number of initially active nodes, which are selected uniformly at random in each community.
\end{abstract}



	\end{frontmatter}

\section{Introduction}

Bootstrap percolation on a graph is a simple activation process
that starts with a given number of initially active nodes (called seeds) and evolves as follows.
Every inactive node that has at least $r\geq 2$ active neighbors is activated, and remains so forever.
The process stops when no more nodes can be activated. There are two main cases of interest: one in which the seeds are selected uniformly at random
among the nodes, and one in which the seeds are arbitrarily chosen. In both cases, the main question concerns the final size of the set of active nodes.
{
Bootstrap percolation was introduced in \cite{chalupa} on a Bethe lattice, and successively investigated
on regular grids and trees \cite{bollobas,BPP}. More recently, bootstrap percolation has been studied
on random graphs and random trees \cite{fountolakis2,amini1,amini2,amini3,angel,pittel,galton,rgg,fountolakis,JLTV,kozma16,turova15}, motivated by the increasing interest in large-scale complex systems such as technological, biological and social networks. For example, in the case of social networks,
bootstrap percolation may serve as a primitive model for the spread of ideas, rumors and
trends among individuals. Indeed, in this context one can assume that a person will adopt
an idea after receiving sufficient influence by friends who have already adopted it \cite{kempe,munik,watts}.

In more detail, bootstrap percolation has been studied on random regular graphs \cite{pittel},
on random graphs with given vertex degrees \cite{amini1}, on Galton--Watson random trees \cite{galton}, on random geometric graphs \cite{rgg}, on Chung--Lu random graphs \cite{amini2,amini3,fountolakis} (\textcolor{black}{which notably permit considering
the case of power-law node degree distribution}),
on small-world random graphs \cite{kozma16,turova15} and on Barabasi--Albert random graphs \cite{fountolakis2}.
Particularly relevant to our work is the paper by Janson et al. \cite{JLTV}, where the authors have provided a detailed analysis
of the bootstrap percolation process on the Erd\H{o}s--R\'{e}nyi random graph. We emphasize that in \cite{JLTV} the seeds are chosen
 uniformly at random among the nodes, however, as proved in \cite{feige}, the critical number of seeds 
 triggering percolation can be significantly reduced if the selection of seeds is optimized.

Over the years, several variants of the bootstrap percolation have 
been considered. In majority bootstrap percolation,
a node becomes active if at least half of its neighbors are active.
In jigsaw percolation, introduced in \cite{jigsaw}, there are two types of edges, one
representing  \lq\lq social links" and one representing
\lq\lq  compatibility of ideas". Two clusters of nodes merge together
if there exists at least one edge of each type between them.
Majority and jigsaw bootstrap percolation have been
analyzed on the Erd\H{o}s--R\'{e}nyi random graph in \cite{cecilia17} and
\cite{belajigsaw}, respectively.

Community structure is an important
characteristic of many real-world graphs. This feature, however, is not
present in any of the graphs on which bootstrap percolation (or its
variants) have been studied so far.
Informally, one says that a graph has a community structure
if nodes are partitioned into clusters in such a way that many edges join nodes of the same cluster and comparatively
fewer edges join nodes of different clusters \cite{girvanpnas02}. Many methods have been proposed for community detection in real networks
(see the review article \cite{fortunato}). 

Through the development of the theoretical
foundations of community detection, the so-called stochastic block
model (SBM) has arisen naturally, and attracted considerable attention.
The SBM is essentially the superposition of Erd\H{o}s--R\'{e}nyi
graphs, and is perhaps the simplest interesting case of a random
graph with community structure.
In particular, detection of two symmetric communities has been studied in \cite{mass14}, while partial or exact recovery of the community
membership has been investigated in \cite{abbe15,abbe16}.
}

In this paper we study classical bootstrap percolation on the SBM with two \textcolor{black}{(in general asymmetric)} communities, assuming
that seeds are selected uniformly at random within each community and allowing
a different number of seeds for different communities. We prove the existence of a sharp phase transition for the number
of eventually active nodes,
identifying a sub-critical regime, in which the evolution of the  bootstrap percolation process is very limited
(in the sense that the final size of active nodes is of the same order as the number of seeds), 
and a super-critical regime, in which
the activation process percolates almost completely
(in the sense that the  vast majority of nodes will be activated).
Although our results generalize some of the main achievements in \cite{JLTV},
we emphasize that our techniques significantly differ from those employed in \cite{JLTV}.
In particular, we devise a suitable extension of the classical {\em  binomial chain} construction
originally proposed in \cite{Scalia} (and also used in \cite{JLTV}), 
adapting it to the SBM. Furthermore, as opposed to \cite{JLTV}, where Doob's martingale inequality is employed,
we use deviation inequalities for the binomial distribution to prove that bootstrap percolation on the SBM concentrates around its average.
Our approach provides exponential bounds on the related tail probabilities,
which allow us
to strengthen the convergence in probability for the final size
of active nodes (as obtained in \cite{JLTV}) to the level of almost sure convergence.

To better understand the main difficulties in the analysis of the bootstrap percolation process on the SBM, we recall
that in the classical binomial chain construction a (virtual) discrete time is introduced: at each time step
a single active node is {\em explored}  by revealing its neighbors. Nodes become active as soon as the number of their  explored neighbors reaches the percolation threshold $r$.
In the SBM the stochastic properties of the set of active nodes at time step $t$ heavily depend on
the number of nodes that have been explored in each community up to time $t$, and this makes the analysis
of the bootstrap percolation process on the SBM
significantly more complex.
In particular, it
requires the identification of an appropriate {\em strategy} to select the community in which a new node is  explored at every time step.

Although considerably flexible and mathematically tractable, the SBM
does not accurately describe most real-world networks.
For instance, it does not allow for heterogeneity of nodes within communities.
Different variants of the SBM have been proposed to better fit real network data, such as letting nodes
follow a given degree sequence \cite{amin10,newman11} or considering overlapping communities and mixed membership models \cite{airoldi08,gopalan13}.
We acknowledge that analyzing bootstrap percolation on the SBM is only a first step
towards a better understanding of this process on more sophisticated
community-based models.

The paper is organized as follows. {\color{black} In Section \ref{sec:SBM} we introduce the model
and our assumptions on its parameters. 
The main results of the paper are stated in Section \ref{sec:main}, together
with some numerical illustrations.
In Section \ref{sec:method} we provide an overview of our analysis, by first
introducing the extension to the SBM of the classical binomial chain representation of the bootstrap percolation process,
and then by giving a high-level description of our proofs.
The detailed proofs are reported in Section \ref{sec:proofs}. 
 {Lastly, in  \textcolor{black}{Appendix}, we report  the proof of some ancillary results.}

\section{The stochastic block model}\label{sec:SBM}

\subsection{Model description}
The SBM $G=G(n_1,n_2,p_1,p_2,q)$,
with  number of nodes $n=n_1+n_2$ and  parameters $p_1,p_2,q\in [0,1)$, 
is a random graph formed by the union of two disjoint Erd\H{o}s--R\'enyi random graphs $G_i=G(\osn{i},p_i)$, $i=1,2$,
called hereafter communities,  where edges joining  nodes in different communities $G_{1}$ and $G_2$ are independently
added with probability $q$.  In the following we will refer to edges between nodes in the same community as 
\lq\lq intra-community" edges and to edges joining nodes in different communities as \lq\lq inter-community" edges.

Bootstrap percolation on the SBM  is an activation process that obeys {to} the following rules:
\begin{itemize}
\item At the beginning, an arbitrary number $a_i$ ($ a_i \leq \osn{i}$) of nodes, called seeds, are chosen
uniformly at random among the nodes of $G_i$. Seeds are declared to be active, while  nodes not belonging to the set of seeds are initially inactive.
\item An inactive node becomes active as soon as at least $r\geq 2$ of its neighbors are active, and then remain active forever,
so that the set of active nodes grows monotonically.
\item The process stops when no more nodes can be activated.
\end{itemize}

The bootstrap percolation process naturally evolves through generations 
of nodes that are sequentially activated.  The initial generation $\cG_0$ is the set of seeds; the first generation $\cG_1$ is composed  
 by all those nodes that are neighbors of at least $r$ seeds;
the second generation $\cG_2$ is composed by all the nodes that are neighbors of at  least $r$ nodes in $\cG_0\cup \cG_1$, and so on. The bootstrap percolation process stops when either
an empty generation is obtained or all the nodes are active. The final set of active nodes is clearly given by
\[
\cG\equiv\bigcup_{k\ge 0 
}	\cG_k. 
\]

We conclude this subsection introducing some notation and terminology.
Given two functions $f_1$ and $f_2$ we write $f_1(m)\ll f_2(m)$ (or equivalently $f_1(m)=o(f_2(m))$), $f_1(m)\sim f_2(m)$,
and $f_1(m)=O(f_2(m))$ if, as $m\to\infty$, $f_1(m)/f_2(m)\to 0$,
$f_1(m)/f_2(m)\to 1$ and $\limsup_{m\to\infty}|f_1(m)/f_2(m)|<\infty$. Letting
$|\mathcal X|$ denote the cardinality of a set $\mathcal X$, we say that the bootstrap percolation process {\em percolates}
whenever $|\cG|=n-o(n)$, that is,  whenever almost all the nodes are  activated.

\subsection{Model assumptions}\label{model-assumptions}

In the following we consider a sequence of SBMs with a growing number of nodes $n$.
We warn the reader that, unless explicitly written, all the limits in this paper are taken as $n\to\infty$.

We assume that the communities $G_1$ and $G_2$ have sizes that  are asymptotically of the same order, i.e.,
\begin{equation}\label{eq:a1a3}
n_1\sim\nu n_2,\quad\text{for some $\nu\in\mathbb{R}_{+}:=(0,\infty)$,}
\end{equation}
and that the \textcolor{black}{inter-community} and the intra-community edge probabilities are asymptotically
of the same order too,   i.e.,
\begin{equation}\label{eq:gratiog}
\text{$q\sim\gamma p_1$,\quad for some $\gamma\in\mathbb{R}_+$,\quad $p_1\sim\mu p_2$,\quad for some $\mu\in\mathbb{R}_+$.}
\end{equation}
Note that since $\gamma>0$ the communities are never isolated.
{Similarly to \cite{JLTV}, we assume  
\begin{equation}\label{eq:hyp2bis}
1/n_i\ll p_i\ll 1/(n_i^{1/r}),\quad i=1,2
\end{equation}
and we define the critical number of seeds,
in correspondence of which the bootstrap percolation process exhibits a phase transition in the Erd\H{o}s--R\'enyi random graph $G(n_i,p_i)$,
by
\begin{equation*}
g_i:=\left(1-\frac{1}{r}\right)\left(\frac{(r-1)!}{n_i p_i^r}\right)^{\frac{1}{r-1}},\quad i=1,2.
\end{equation*}
}
As proved in \cite{JLTV}, under
 \eqref{eq:hyp2bis}, we have
 \begin{equation}\label{eq:ptcto0bm}
 g_i\to\infty,\quad g_i/n_i\to 0,\quad p_i g_i\to 0,\quad i=1,2.
 \end{equation}
Note that by \eqref{eq:a1a3} and \eqref{eq:gratiog} it follows that $g_1$ and $g_2$ are asymptotically comparable.
Furthermore, similarly to \cite{JLTV},  we assume
{
\begin{equation}\label{eq:trivial}
a_i/g_i \to \alpha_i\geq 0,\quad i=1,2, \quad\text{with $\max\{\alpha_1,\alpha_2\}>0$.}
\end{equation}
Without loss of generality, we suppose
\begin{equation}\label{eq:alfaorder}
\alpha_1\ge\alpha_2\quad\text{with $\alpha_1>0$.}
\end{equation}
}

Inspired by some literature on the subject (see e.g. \cite{Newman}) we say that the SBM is {\em  assortative}
{if the intra-community edge probabilities exceed the inter-community edge probability.
Specifically, a SBM is said assortative if $q^2<p_1 p_2$.
Since  $q^2/(p_1p_2)\to\gamma^{2}\mu$ (see \eqref{eq:gratiog}),  by setting 
\[
\chi_{ii}=1,\quad i=1,2,\quad
\chi_{12}:=\gamma(\nu\mu^r)^{\frac{1}{r-1}}, \quad\chi_{21}:=\gamma(\nu\mu)^{-1/(r-1)}\quad\text{and}\quad\bm{\chi}=(\chi_{ij})_{i,j=1,2},
\]
the {\em assortative} condition can be reformulated as
$\mathrm{det}\bm{\chi}>0$. }Therefore, in the following we will refer to  {\it assortative} {SBM} when  $\mathrm{det}\bm{\chi}>0$, {\it dis-assortative} {SBM} when $\mathrm{det}\bm{\chi}<0$ and
 {\it neutral} {SBM} when $\mathrm{det}\bm{\chi}=0$. {Although these notions do not play any role in our main  results
(i.e., Theorems \ref{teo-blocksub} and \ref{teosupcrit}), they do 
have an impact on the definition of the {\it critical curve} for the system (see Proposition \ref{prop:alpharegion}).}

Finally, we remark once again that, within each community, the seeds  must be selected uniformly at random and in such a way that
the number of seeds satisfies the constraint \eqref{eq:trivial}.

\subsection{Bootstrap percolation on the Erd\H{o}s--R\'enyi random graph: a quick review}

To better position our results with respect to the existing literature, we briefly recall the main {achievements} 
in \cite{JLTV}. Note that the Erd\H{o}s--R\'enyi random graph corresponds to a SBM with a single community, (i.e.,
$i=1$). It has been proved in \cite{JLTV} (see Theorem 3.1$(ii)$). that:\\
\noindent $(i)$ If \eqref{eq:hyp2bis} and \eqref{eq:trivial} hold (with $i=1$) and $\alpha_1<1$, then
\[
|\cG|/g_1\to\frac{r\varphi(\alpha_1)}{(r-1)\alpha_1},\quad\text{in probability}
\]
where $\varphi(\alpha_1)$ is the unique solution in $[0,1]$ of { equation} $rx-x^r=(r-1)\alpha_1$ with unknown $x$ (see Theorem 3.1$(i)$ in \cite{JLTV}).\\
\noindent $(ii)$ If \eqref{eq:hyp2bis} and \eqref{eq:trivial} hold (with $i=1$) and $\alpha_1>1$, then
\[
|\cG|/n\to 1,\quad\text{in probability}.
\]

\section{Main results}\label{sec:main}

The bootstrap percolation process on the Erd\H{o}s--R\'enyi random graph exhibits a sharp phase transition, see \cite{JLTV}.
The reader may be wondering whether more complex phenomena, such as selective percolation of communities,
can be observed on the SBM. We will show that this is not the case. Indeed, under the assumptions described in Subsection \ref{model-assumptions},
the bootstrap percolation process either stops with high probability when $O(g_1)$ vertices have been activated (sub-critical case) or  percolates (super-critical case).  
A selective percolation of the communities may be instead  observed when $\gamma=0$ (i.e.,  $q=o(p_1)$), where the bootstrap percolation process may behave
in each community as if they were isolated.
%

To state our main results we need some additional notation. {For $\mathbf{x}=(x_1,x_2)\in [0,\infty)^2$,}
we define the following functions:
\[
\rho_i(\mathbf{x}):=\alpha_i-x_i+r^{-1}(1-r^{-1})^{r-1}(x_i+\chi_{ij}x_j)^r,\quad i \neq j \in\{1,2\}  
\]
and the following sets: 
\[
\mathcal{D}:=\left\{\mathbf{x}\in [0,r/(r-1)]^2:\,\,x_1+\chi_{12}x_2\leq\frac{r}{r-1},\,\,x_2+\chi_{21}x_1\leq\frac{r}{r-1}\right\},
\]
\[
\oE_1:=\{\mathbf{x}\in\oD:\,\,\rho_1(\mathbf{x})\leq 0\},\quad\oE_2:=\{\mathbf{x}\in\oD:\,\,\rho_2(\mathbf{x})\leq 0\},
\]
\[
\widetilde\oE_1 :=\{\mathbf{x}\in\oD:\,\,\rho_1(\mathbf{x})=0\},\quad\widetilde \oE_2:=\{\mathbf{x}\in \oD:\,\,\rho_2(\mathbf{x})=0\}.
\]
For a set $\mathcal{H}\subset\mathcal{D}$, we denote by $\overset\circ{\mathcal{H}}$ its interior (with respect to the Euclidean topology on $\mathbb{R}^{2}$ restricted
to $\mathcal{D}$). Throughout this paper, we consider the following three disjoint and exhaustive conditions:\\
\noindent$(\bold{Sub})$: $\overset\circ{\mathcal{E}}_1\cap\overset\circ{\mathcal{E}}_2\neq\emptyset$,
\qquad$(\bold{Crit})$: $\overset\circ{\mathcal{E}}_1\cap\overset\circ{\mathcal{E}}_2=\emptyset$,
$\widetilde{\mathcal{E}}_1\cap\widetilde{\mathcal{E}}_2\neq\emptyset$,\qquad$(\bold{Sup})$: $\mathcal{E}_1\cap\mathcal{E}_2=\emptyset$.\\
Hereafter, we refer to such conditions as sub-critical, critical and super-critical regimes, respectively.
A graphical representation of these regimes  is given in Figures \ref{fig:fig_sub}, \ref{fig:fig_crit} and \ref{fig:fig_sup}, where the blue curves represent $\widetilde{\mathcal{E}}_1$ and $\widetilde{\mathcal{E}}_2$ (additional notation appearing on the plots will be introduced later on).

\begin{remark}\label{re:alfamag1}
Let $i\in\{1,2\}$ be fixed. A straightforward computation shows that if $\alpha_i>1$, then $\min_{\bold x\in\mathcal{D}}\rho_i(\bold x)>0$, therefore $\mathcal{E}_i=\emptyset$ and so $(\bold{Sup})$ holds. \textcolor{black}{Consequently, conditions $(\bold{Sub})$ and $(\bold{Crit})$  imply 
 $\alpha_i\leq 1$ for any $i=1,2$.} 
\end{remark}

\subsection{Phase transition on the SBM model}
Next theorems provide the main results of the paper. Hereon, for ease of notation, we denote by $(C)$ the set of conditions: \eqref{eq:a1a3},  \eqref{eq:gratiog},
\eqref{eq:hyp2bis}, \eqref{eq:trivial} and \eqref{eq:alfaorder}.

\threplots{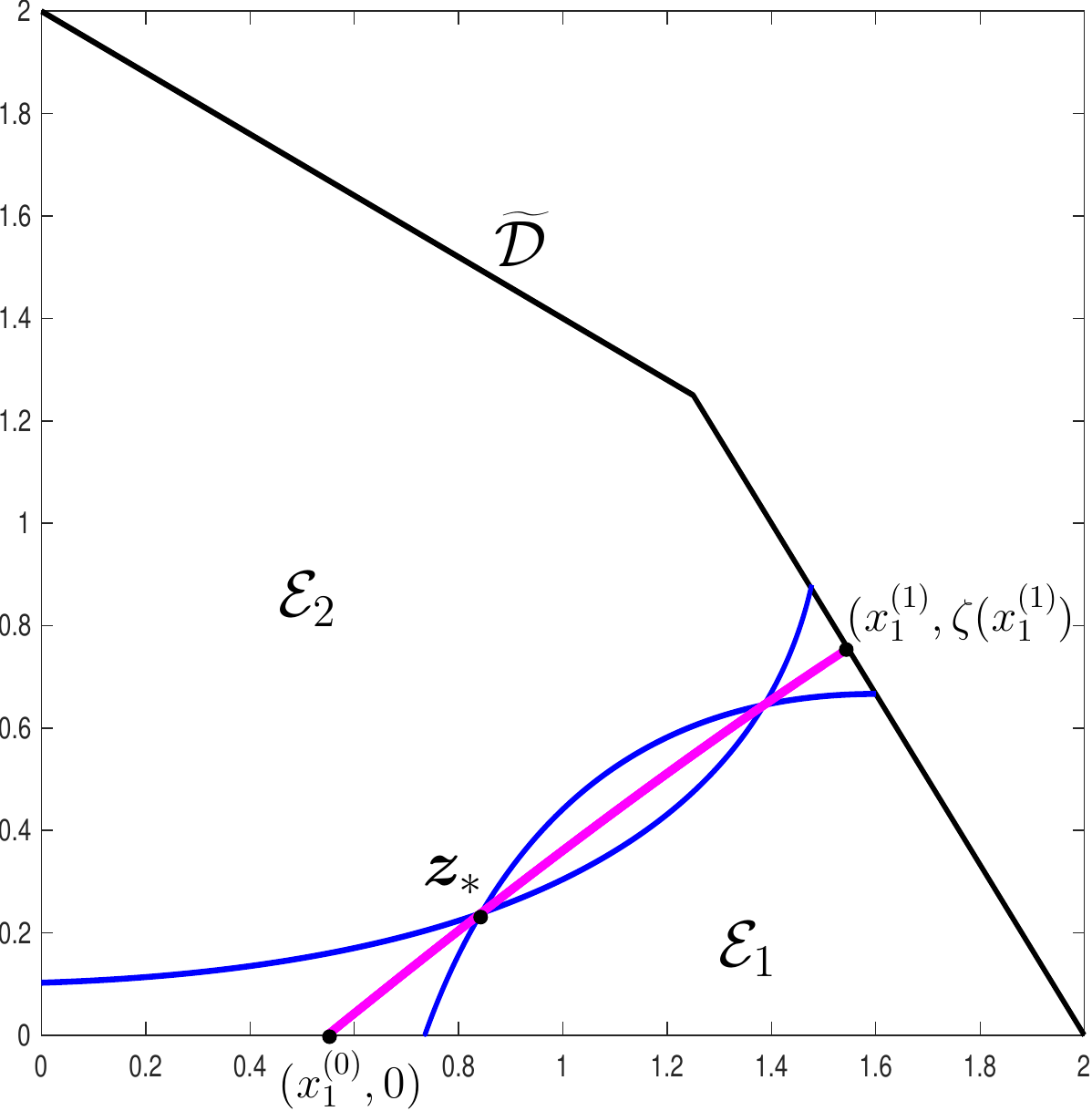}{Sub-critical regime,  \mbox{$r=2$}, $\chi_{12}=\chi_{21}=0.6$,
\mbox{$\alpha_1=0.56$}, $\alpha_2=0.1$.}{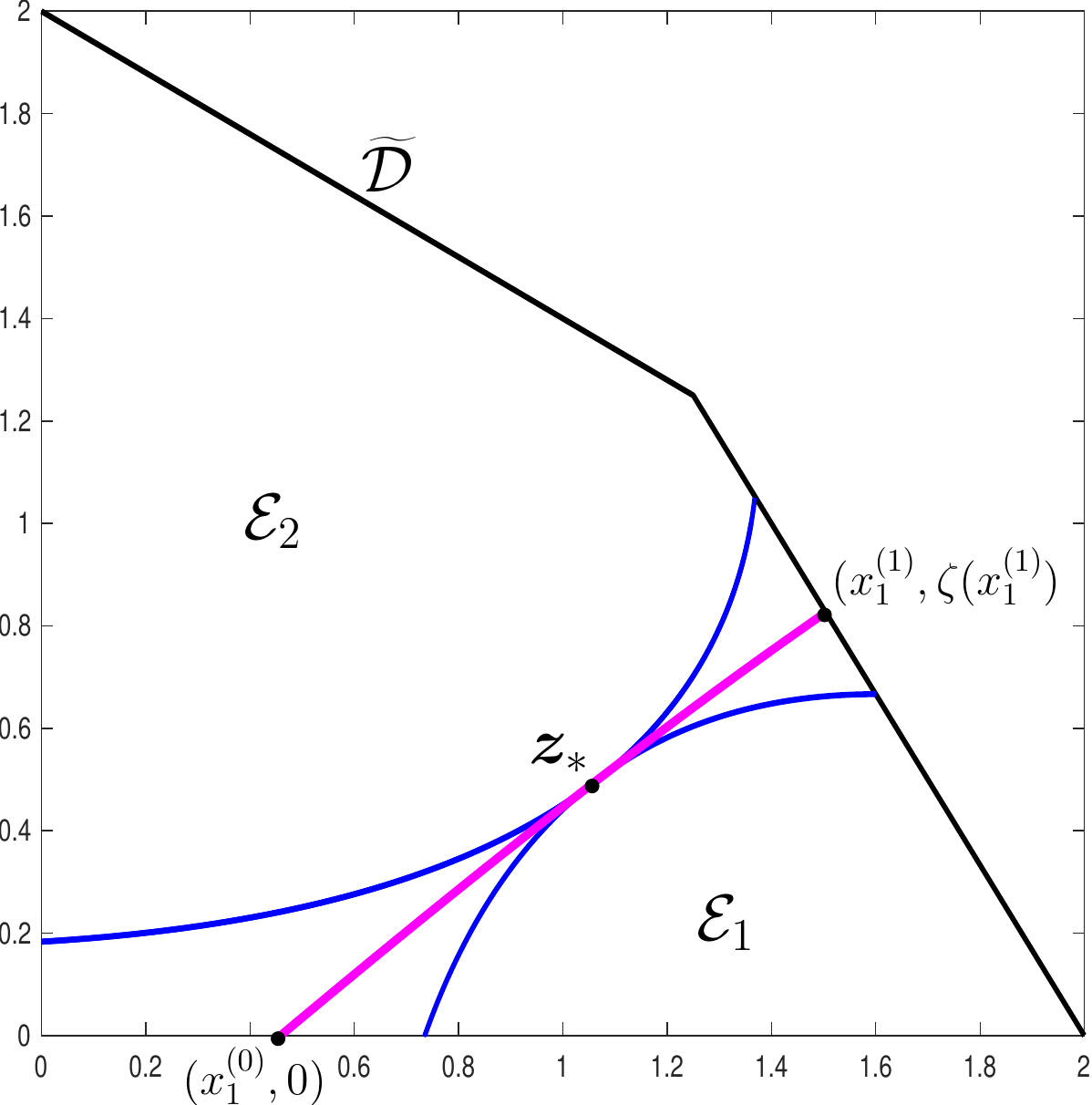}{Critical regime,  \mbox{$r=2$}, $\chi_{12}=\chi_{21}=0.6$, \mbox{$\alpha_1=0.6$}, $\alpha_2=0.175$. }{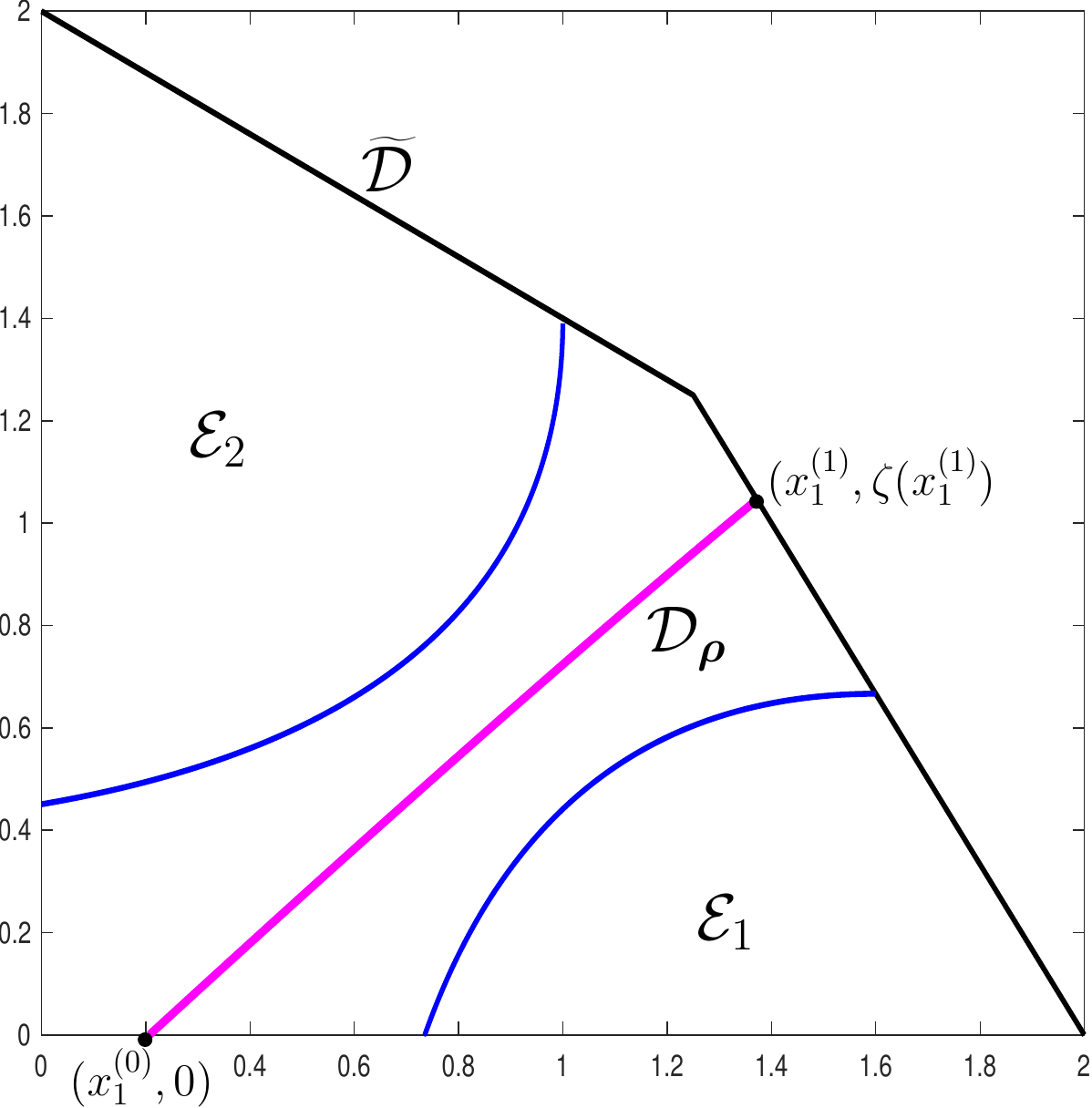}{Super-critical regime,  $r=2$, $\chi_{12}=\chi_{21}=0.6$, \mbox{$\alpha_1=0.6$}, $\alpha_2=0.4$.}

\begin{theorem}\label{teo-blocksub}
Assume $(C)$ and $(\bold{Sub})$. Then, for any $\varepsilon>0$ there exists $c(\varepsilon)\in\mathbb R_+$ such that
\begin{equation}\label{eq:limsub}
P\left(\Big|\frac{|\mathcal{G}|}{g_1}-x_*\Big|>\varepsilon\right)=O(\mathrm{e}^{-c(\varepsilon)g_1}),
\end{equation}
where the explicit expression of  the positive constant $x_*>0$ is given in \eqref{eq:xstar}.
\end{theorem}
\begin{theorem}\label{teosupcrit}
Assume $(C)$ and $(\bold{Sup})$. Then, for any $\varepsilon>0$ there exists $c(\varepsilon)\in\mathbb R_+$ such that
\begin{equation}\label{eq:limsup}
P\left(\Big|\frac{|\mathcal{G}| }{n}-1\Big|>\varepsilon\right)=O(\mathrm{e}^{-c(\varepsilon)g_1}).
\end{equation}
\end{theorem}

Roughly speaking, the above results can be rephrased as follows:\\
\noindent $(i)$ under $(C)$ and $(\bold{Sub})$, the bootstrap percolation process on the SBM reaches,
as $n\to\infty$, a final size of active nodes which is of the same order as $a_1+a_2$
(indeed, by \eqref{eq:gratiog}, the definition of $g_i$,
\eqref{eq:trivial} and \eqref{eq:alfaorder}, it easily follows that
$a_1+a_2\sim (\alpha_1+\alpha_2(\nu\mu^r)^{1/(r-1)})g_1$),\\
\noindent $(ii)$ under $(C)$ and $(\bold{Sup})$, the bootstrap percolation process on 
the SBM percolates, as $n\to\infty$.

\begin{remark}
Replacing the assumption \eqref{eq:hyp2bis} with the (slightly) stronger condition: 
\begin{equation*}
\text{For any $i=1,2$, $1/n_i\ll p_i$ and either $p_i\ll 1/(n_i)^\frac{1}{r'}$ or
$p_i\sim c/(n_i)^\frac{1}{r'}$, for some $c>0$ and $r'\in (r,\infty)$,}
\end{equation*}
\textcolor{black}{by a standard application of the Borel-Cantelli lemma,
the claims \eqref{eq:limsub} and \eqref{eq:limsup}
can be strengthened, respectively,  as
$\frac{|\mathcal{G}|}{g_1} \to x_*$ and $\frac{|\mathcal{G}|}{n} \to 1$ almost surely.} 
\end{remark}

\subsection{The critical curve and the sub-critical and super-critical regions}
To complement the results of Theorems \ref{teo-blocksub} and \ref{teosupcrit}, in this subsection we determine the sub-critical and the super-critical regions of the system, i.e.,  
the set of $\bm{\alpha}=(\alpha_1,\alpha_2)$ for which either the sub-critical or the super-critical behavior is observed. 
We restrict our investigation to $\bm{\alpha}\in[0,1]^2$ since, as already noticed in Remark \ref{re:alfamag1},
\textcolor{black}{ necessarily $(\bold{Sup})$  holds whenever  $\alpha_1>1$ and/or $\alpha_2>1$.}  We write $\overset\circ{\mathcal{E}}_{i}(\alpha_i)$, $\widetilde{\oE}_i(\alpha_i)$ and $\oE_i(\alpha_i)$
in place of $\overset\circ{\mathcal{E}}_{i}$, $\widetilde{\oE}_i$ and $\oE_i$, 
respectively,  to  make the dependence on $\alpha_i$ explicit. We define the regions
\begin{equation}
\mathcal{R}_{\bold{Sub}}:=\left\{\bm{\alpha}\in [0,1]^2:\,\,\overset\circ{\mathcal{E}}_1(\alpha_1)\cap\overset\circ{\mathcal{E}}_2(\alpha_2)\neq\emptyset \right\},\quad
\mathcal{R}_{\bold{Sup}}:=\left\{\bm{\alpha}\in [0,1]^2:\,\, {\mathcal{E}}_1(\alpha_1)\cap{\mathcal{E}}_2(\alpha_2)=\emptyset  \right\},\label{Rcond-sup}
\end{equation}
and the curve
\begin{equation}
\mathcal{R}_{\bold{Crit}}:=\left\{\bm{\alpha}\in [0,1]^2:\,\, \overset\circ{\mathcal{E}}_1(\alpha_1)\cap\overset\circ{\mathcal{E}}_2(\alpha_2)=\emptyset,
\widetilde{\mathcal{E}}_1(\alpha_1)\cap\widetilde{\mathcal{E}}_2(\alpha_2)\neq\emptyset \right\}, \label{Rcond-sub-crit}
\end{equation}
to which we refer as the sub-critical and the super-critical regions, and the critical curve, respectively.

By exploiting the convexity of the functions $\rho_i(\cdot)$ and by imposing the tangency condition between the curves $\widetilde{\mathcal E}_1(\alpha_1)$ and $\widetilde{\mathcal E}_2(\alpha_2)$, one can show 
the following Proposition \ref{prop:alpharegion}, whose proof is elementary, and therefore omitted.  From here on, we denote by $\mathcal M^t$ the transpose of the matrix $\mathcal M$.

\begin{proposition}\label{prop:alpharegion}
{The following claims hold}:\\
\noindent$(i)$ Under $(C)$ and $\mathrm{det}\bm{\chi}\neq 0$, we have 
\begin{align*}
\mathcal{R}_{\bold{Crit}}&=\Biggl\{(y_1,y_2)(\bm{\chi}^{-1})^t -r^{-1}(1-r^{-1})^{r-1}(y_1^{r},y_2^r)\in [0,1]^{2}:\,\,0\leq y_1\leq r/(r-1),\nonumber\\
&\qquad\qquad\qquad
y_2=(1-r^{-1} )\left(\frac{1-(1-r^{-1})^{r-1}  y_1^{r-1}}{1 - (1-r^{-1})^{r-1}y_1^{r-1}\mathrm{det}\bm{\chi}}\right)^{1/(r-1)}\Biggr\}.\nonumber
\end{align*}
\noindent$(ii)$ Under $(C)$ and $\text{det}\bm{\chi}=0$,  we have 
\begin{align*}
\mathcal{R}_{\bold{Crit}}&=\Biggl\{(y_1,y_2)-r^{-1}(1-r^{-1})^{r-1}((y_1+y_2\chi_{12})^r,(y_1\chi_{21}+y_2)^r)\in [0,1]^2:\,\,0\leq y_1\leq r/(r-1),\nonumber\\
&\qquad\qquad\qquad  
y_2=\chi_{21}\left[\frac{r}{r-1}\left(\frac{1}{1+\chi_{21}^{r-1}}\right)^{1/(r-1)}-y_1\right] \Biggr\}.
\end{align*}
\noindent$(iii)$ Under $(C)$, $\mathcal{R}_{\bold{Sub}}$ is the convex set delimited
by the curve $\mathcal{R}_{\bold{Crit}}$ and the coordinate axes. 
\end{proposition}

Note that $\mathcal{R}_{\bold{Crit}}$ depends only on the asymptotic properties of the SBM, which are expressed in terms 
of the parameters $r$, $\gamma$, $\mu$ and $\nu$.  In other words, two (sequences of) SBMs with the same parameters {$r$, $\gamma$, $\mu$ and $\nu$} lead to the same critical curve 
$\mathcal{R}_{\bold{Crit}}$,  {and therefore to the same sub-critical and super-critical regions}. 
Hereafter, we illustrate numerically Proposition \ref{prop:alpharegion}, taking  (sequences of) SBMs with parameters $r=2$, $\gamma=0.25$, and $\nu=\mu=1$
 as baseline case.

We start by investigating the impact  of the various parameters  on the sub-critical and the super-critical regions.
To this aim, we vary a parameter at a time, keeping fixed all the others, and determine the critical curve.   

\sidebyside{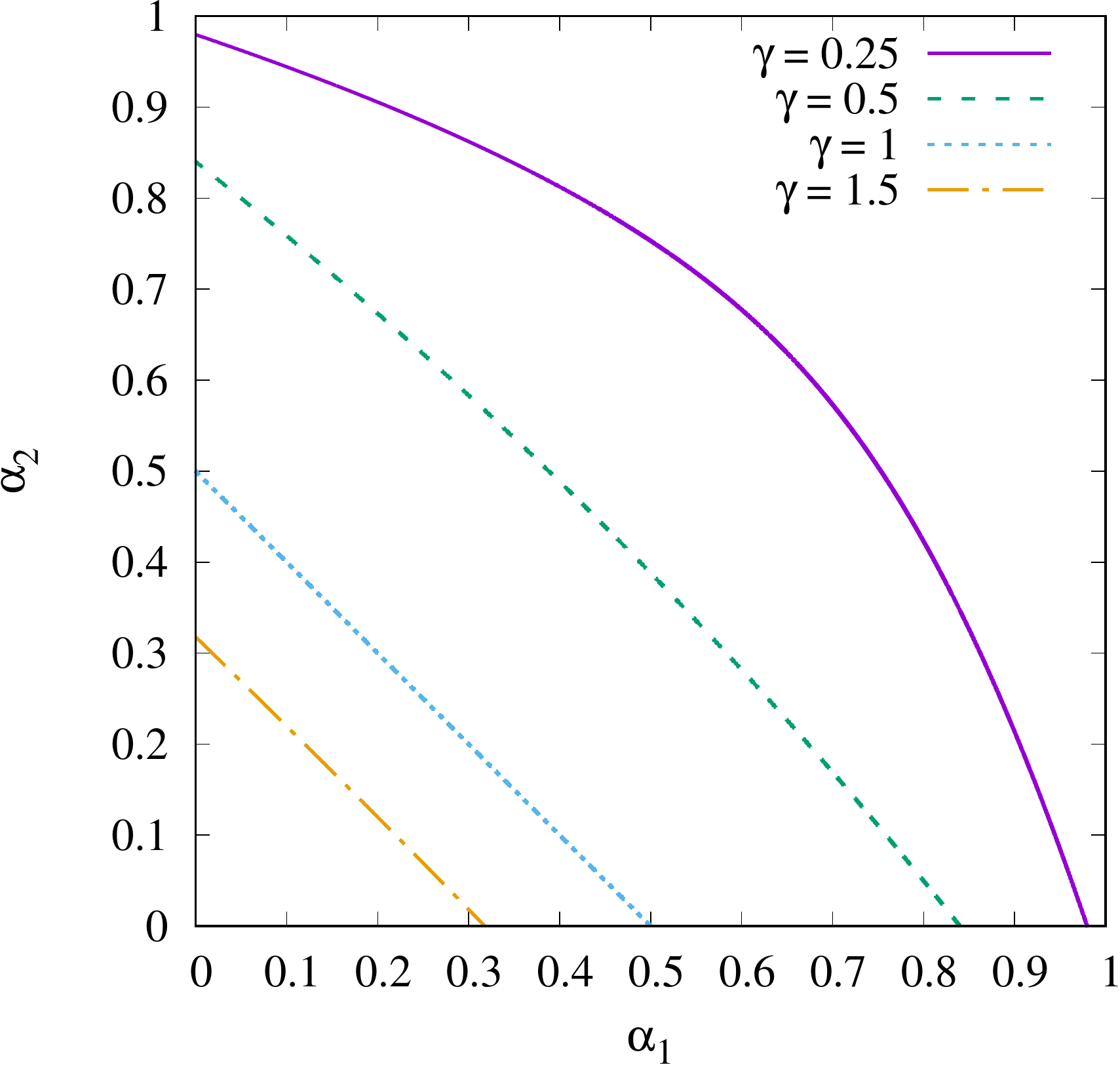}{Critical curves $\mathcal{R}_{\bold{Crit}}$ for different values of $\gamma$.}
{rvary}{Critical curves $\mathcal{R}_{\bold{Crit}}$ for different values of $r$.}

In Figure \ref{fig:gammavary} we vary the parameter $\gamma$, which 
characterizes the strength of the inter-community connectivity with respect to the intra-community connectivity.
When $\gamma<1$   ({$r=2$}, $\mu=\nu=1$) SBMs are assortative, whereas when $\gamma>1$ SBMs are disassortative. 
Finally, in the special case when $\gamma=1$ the SBMs are neutral (i.e., $\text{det} {\bm \chi}= 0$)
and  exhibit the same $\mathcal{R}_{\bold{Crit}}$ of Erd\H{o}s--R\'{e}nyi random graphs. 
In this special case the critical curve corresponds to the segment
where $\alpha_1 + \alpha_2 = 0.5$,  {indeed a straightforward computation gives
$g_1=\frac{1}{2 n_1 p_1^2}$}.
We further note that, as $\gamma\downarrow 0$, the sub-critical region approaches the whole square 
{(because,  as $\gamma\to 0$, the fraction of edges connecting the two communities tends to vanish, 
and therefore the activation process spreads in  the two communities as if they were isolated). }
\textcolor{black}{	Finally, since the sub-critical region is convex for any $\gamma$, in a
 SBM with $\mu=\nu=1$ (i.e., symmetric), we have that the critical number of seeds 
 is minimized when all the seeds are placed in the same community 
(i.e., either $\alpha_1=0$ or $\alpha_2=0$). Instead, the critical number of seeds is maximized 
 when the seeds are equally partitioned  between the communities 
(which approximately occurs, notably, when the seeds are chosen uniformly 
at random among the nodes).} Interestingly, in the latter case (i.e., when the seeds are equally partitioned between the communities), 
a simple computation shows that the critical threshold 
in a SBM with $\mu=\nu=1$ is asymptotically equal to the critical threshold in an Erd\H{o}s--R\'{e}nyi  random graph having the same average degree.

In Figure \ref{fig:rvary} we vary the threshold parameter $r$.   
Note that,
as $r\uparrow\infty$, the sub-critical region approaches the domain 
\[
\{(\alpha_1,\alpha_2)\in [0,1]^2:\,\,\alpha_1+\gamma\,\alpha_2<1,\,\,\alpha_2+ \gamma\,\alpha_1<1\}
\]
\textcolor{black}{(this property holds for any $\gamma$ in the symmetric SBM with $\mu=\nu=1$). }

\sidebyside{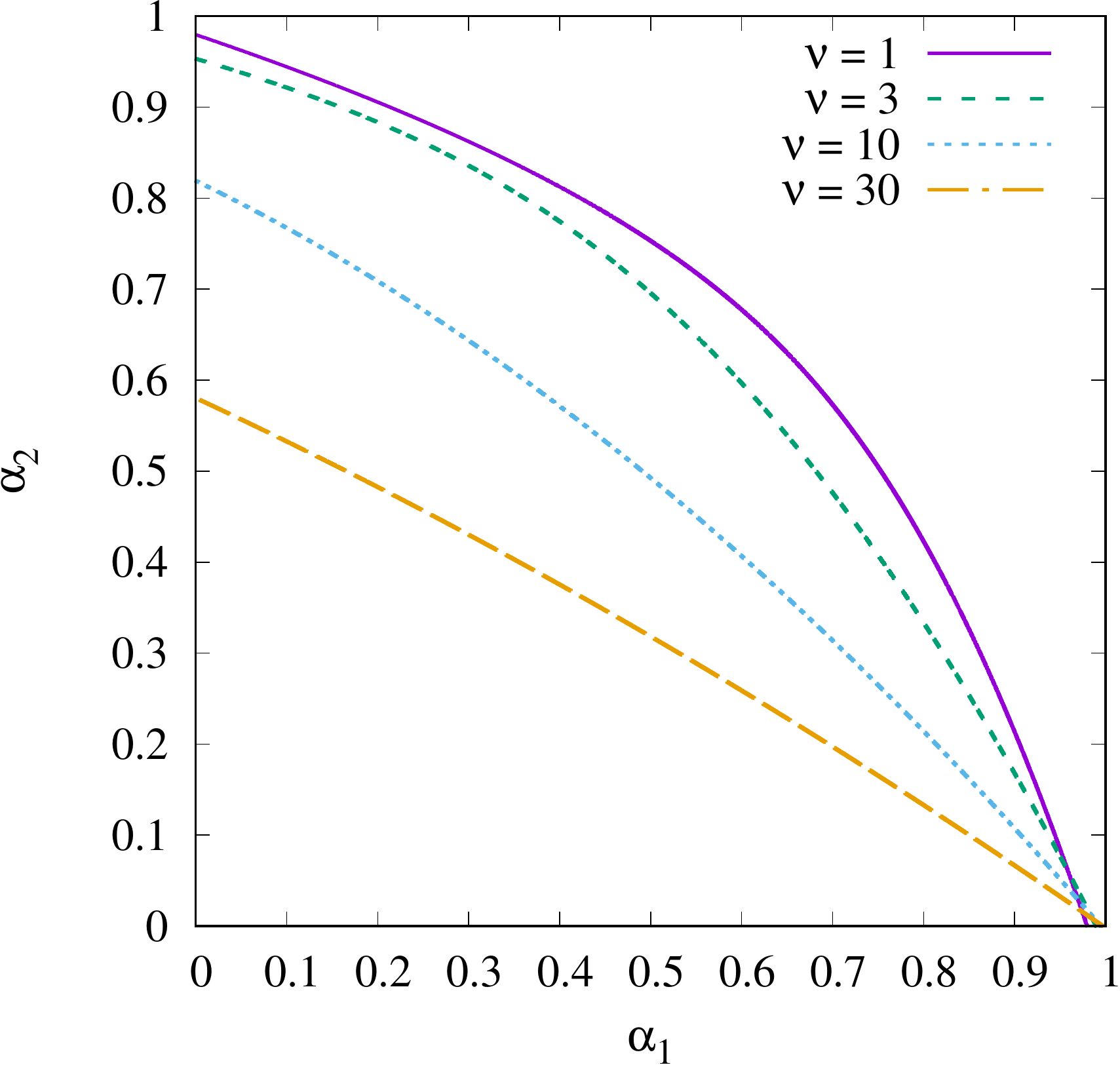}{Critical curves $\mathcal{R}_{\bold{Crit}}$ for different values of $\nu$.}
{muvary}{Critical curves $\mathcal{R}_{\bold{Crit}}$ for different values of $\mu$.}

Next, we explore what happens in SBMs with $\mu\neq\nu$ (i.e., asymmetric) by changing either $\nu$ or $\mu$. 
In Figure \ref{fig:nuvary} we {fix $r=2$, $\gamma=0.25$, $\mu=1$} and increase the parameter $\nu$, making the first community increasingly larger
than the second community.  \textcolor{black}{Interestingly,
we observe a significant reduction in the (normalized) critical value of $\alpha_2$ 
for increasing values of $\nu$ when we put all the seeds in the community $G_2$ (i.e., $\alpha_1 = 0$): this means that fewer and fewer
seeds are needed in community $G_2$ to trigger percolation, as the community $G_1$ becomes larger and larger. This because the epidemic
transfers into the community $G_1$, where it propagates more easily thanks to the larger number of available nodes.
 However, note that, to minimize the  (un-normalized) critical number of seeds,
all the seeds must be placed in the larger community $G_1$,  as a consequence of the fact that $g_i$, $i=1,2$, are different.}

Finally, in Figure \ref{fig:muvary} we {fix $r=2$, $\gamma=0.25$, $\nu=1$} and increase the parameter $\mu$,  thus increasing the 
intra-community probability in $G_1$.
For large values of $\mu$,  considerations similar to  Figure \ref{fig:nuvary} apply.

\section{Overview of the methodology}\label{sec:method}
\subsection{The extension of the binomial chain construction}
We introduce  a  discrete time $t\ge  0$  
and we assign a marks counter $M_{v}(t)$, $M_v(0):=0$, to every node $v$ which is not a seed. Seeds are activated at time $t=0$. 
We set $\mathcal{U}_{i}(0):=\emptyset$ and denote by $\mathcal{A}_{i}(0)$ the set of seeds in the community $G_i$.
The process, then,  evolves according to the following recursive procedure. At   time $t\in\mathbb N:=\{1,2,\ldots\}$:

\begin{itemize}
\item We  arbitrarily  select   a community $G_j$ provided that  $\mathcal{A}_{j}(t-1)\setminus\mathcal{U}_{j}(t-1)\neq\emptyset$.
\item From the selected community $G_j$, we choose, uniformly at random, a node $v\in\mathcal{A}_{j}(t-1)\setminus\mathcal{U}_{j}(t-1)$.
\item We use the chosen node $v$, i.e., we explore the node $v$ by revealing its neighbors and by adding a mark to each of them.
\item We set $\mathcal{U}_{j}(t):=\mathcal{U}_{j}(t-1)\cup\{v\}$ and $\mathcal{U}_{i}(t):=\mathcal{U}_{i}(t-1)$, for $i\neq j$.
We also set $\mathcal{A}_{i}(t):=\mathcal{A}_{i}(t-1)\cup\Delta\mathcal{A}_{i}(t)$,  
where
$\Delta\mathcal{A}_{i}(t)$ is the set of nodes in the community $G_i$ that become active exactly at time $t$, i.e., the set of nodes in $G_i$
that have received the $r$-th mark exactly at time $t$. Note that  $\Delta\mathcal{A}_{i}(t)=\emptyset$ for  $t<r$,
since no other nodes are activated until at least $r$ seeds are used.
\item The process terminates as soon as there are no active and  still  unused nodes, i.e., at time step:
\begin{equation}\label{eq:T}
T:=\min\{t\in\mathbb N:\,\,\mathcal{A}_{i}(t)=\mathcal{U}_{i}(t),\,\forall i=1,2\}.
\end{equation}
\end{itemize}

{Note that, since only one node is used at each time step, for any $t\leq T$,}
$|\mathcal{U}(t)|=t$,  where $\quad\mathcal{U}(t):=\mathcal{U}_{1}(t)\cup\mathcal{U}_{2}(t)$.
Let $\mathcal{A}(t):=\mathcal{A}_1(t)\cup\mathcal{A}_2(t)$ denote the set of active nodes at time $t\leq T$.
{We clearly have
\begin{equation}\label{eq:AMv}
v\in\mathcal{A}(t)\setminus\mathcal{A}(t-1)\quad\text{if and only if}\quad M_v(t)=r,\quad 1\leq t\leq T
\end{equation}
where
\begin{equation}\label{eq:Mv}
M_{v}(t)=\sum_{i=1}^{2}\sum_{s=1}^{U_{i}(t)}I^{(i)}_{v}(s),\quad\text{$\forall$ $v\not\in\mathcal{A}(t-1)$}
\end{equation}}
$U_{i}(t):=|\mathcal{U}_{i}(t)|$ and the random variables $\{I_{v}^{(i)}(s)\}_{v\notin\cA(t),1\leq i\leq 2,1\leq s\leq T}$
are independent, with $I_{v}^{(i)}(s)$ distributed as $\mathrm{Be}(p_{i})$ \footnote{Here $\mathrm{Be}(p)$ denotes a Bernoulli distributed random variable with mean $p\in [0,1]$. }
if $v$ is a node of the community $G_i$, and
distributed as $\mathrm{Be}(q)$ if $v$ is a node of the community $G_j$, $j\neq i$.

{The next proposition guarantees that
the order in which active nodes are used does not have any impact on the 
final set of active vertices $\mathcal{G}$.}
\begin{proposition}\label{prop:equiv}
We have
$
\cG\equiv\mathcal{A}(T).
$
\end{proposition}
{Although Proposition \ref{prop:equiv} may appear rather obvious, it plays a crucial role 
in our proofs. Therefore, for completeness, we report its proof in  Appendix \ref{subsec:prop21}. }

\begin{remark}\label{re:strategy}
{
In the description of the binomial chain representation of the bootstrap percolation process, we did not fully specify the rule according to which  
a community is selected at every time step $t\in\mathbb{N}$. Indeed, we limited ourselves just to mention
a general guideline for the selection of the community: at time $t\in\mathbb{N}$, we choose a community $G_i$ 
which has active and unused nodes. Clearly, this choice can be made in many different  ways.
 Throughout this paper, we refer to such different ways as \lq\lq strategies".
Remarkably, Proposition \ref{prop:equiv} applies to any strategy. It will become clear later on that the opportunity to \lq\lq arbitrarily" 
define the strategy for the community selection, provides a fundamental degree of flexibility that comes in handy  when we analyze the 
bootstrap percolation process on the SBM (see the proofs of Theorems \ref{teo-blocksub} and \ref{teosupcrit}).}
\end{remark}
\textcolor{black}{Hereon, we put $[n]:=\{1,2,\cdots,n\}$} and
let $i\in\{1,2\}$ be fixed. 
We have defined the random marks $I_v^{(i)}(s)$
for $v\notin\cA(t)$ and $1\leq s\leq T$, but, similarly to \cite{JLTV},  see Section 2 therein,
it is possible to introduce additional, redundant random marks, which are independent and Bernoulli distributed with mean $p_i$ if $v$ is a node of the community $G_i$ and with mean $q$ if $v$ is a node of the community $G_j$, $j\neq i$,
in such a way that $I_v^{(i)}(s)$ is defined for all $v\in G$ and \textcolor{black}{$s\in [n]$}.
Such additional random marks are added, for any $1\leq s\le T$,
to already active nodes and so they have no effect on the underlying bootstrap percolation process.
 {This corresponds to artificially extending the chain construction beyond $T$, by selecting and exploring at every time 
 	$T\le t \le n$ a potentially  non-active node.  Clearly such extension has no effect on the dynamics of the bootstrap percolation process up to time $T$, and it is just instrumental.}
Throughout this paper,
we denote by $\mathrm{Bin}(u,p)$, $u\in\mathbb N$, $p\in [0,1]$, a random variable following the binomial distribution with parameters $(u,p)$.

Note that, conditionally  on $U_1(t)$ and $U_2(t)$, the random variable $M_v(t)$ is the sum of independent random variables with the binomial distribution,
i.e.,  for fixed $i\in\{1,2\}$ and 
\textcolor{black}{$t\in [n]\cup\{0\}$} we have
\begin{equation}\label{eq:sumbin}
\text{$M_{v}(t)\,|\,\{\bold{U}(t)=\oboldm(t)\}\overset{\mathcal L}=\mathrm{Bin}(u_i(t),p_i)+\mathrm{Bin}(u_j(t),q)$,\quad $v\in G_i$, $j\neq i$}
\end{equation}
where $\bold{U}(t):=(U_1(t),U_2(t))$, $\oboldm(t):=(u_1(t),u_2(t))$, the symbol $\overset{\mathcal L}=$ denotes the equality in law and the random variables
$\mathrm{Bin}(u_i(t),p_i)$ and $\mathrm{Bin}(u_j(t),q)$ are independent. The number of active nodes in the community $G_i$ at time $t\in [n]\cup\{0\}$
is given by
\begin{equation}\label{eq:Ai}
A_{i}(t):=|\mathcal{A}_{i}(t)|=a_{i}+S_{i}(t),
\end{equation}
where
\begin{equation}\label{eq:Si}
S_{i}(t):=\sum_{v\in G_i\setminus\mathcal{A}_{i}(0)}\bold{1}\{Y_{v}\leq t\},\qquad Y_{v}:=\min\{s\in\mathbb N:\,\,M_{v}(s)\geq r\}.
\end{equation}
Since the random variables
$\{M_{v}(t)\,|\,\{\bold{U}(t)=\oboldm(t)\}\}_{v\in G_i}$
are independent and identically distributed with law specified by \eqref{eq:sumbin}, we have
\begin{equation}\label{eq:bin}
S_{i}(t)\,|\,\{\bold{U}(t)=\oboldm(t)\}\overset{\mathcal L}=\mathrm{Bin}(\osn{i}-a_{i},b_i(\oboldm(t))),
\end{equation}
where
\begin{equation}\label{eq:pi}
b_i(\oboldm(t)):=P\left(\mathrm{Bin}(u_i(t),p_i)+\mathrm{Bin}(u_j(t),q)\geq r\right),\quad i\in\{1,2\},\,j\neq i.
\end{equation}
Hereafter, we denote by $A(t):=|\mathcal{A}(t)|=\sum_{i=1}^{2}A_i(t)$, the number of active nodes in the SBM at time $t$. Note that $|\mathcal{G}|=A(T)=T-1$.

\begin{remark}\label{rem1}

{	
The analysis of the bootstrap percolation process is significantly more complex on the SBM than on the Erd\H{o}s--R\'enyi random graph. Indeed, on
the SBM, 
for any  $t<T$,  the random variables $\{A_i(t)\}_{1\leq i\leq 2}$ depend on the quantities $\{U_i(t)\}_{1\leq i\leq 2}$,  and so on the chosen strategy.
In turn, the choice of a strategy 
is constrained by the availability of active and unused nodes in the different communities.
As a result,
$S_i(t)$ is binomial only given the event $\{\bold{U}(t)=\oboldm(t)\}$. 
In contrast, on the Erd\H{o}s--R\'enyi random graph
the number of used nodes at time $t$ is equal to $t$, and therefore
the law of the number of active and unused nodes at time $t$ is (unconditionally) binomial. 
}
\end{remark}

\subsection{High level description of the proofs} \label{sec:high-level}

{\color{black}
In broad terms, the proofs of Theorems \ref{teo-blocksub} and \ref{teosupcrit} adopt the following approach. First,  note that since \mbox{$T-1=|\mathcal{G}|$}, we can reduce the computation of the tail probabilities of $|\mathcal G|$ to the computation 
of the tail probabilities of $T$.  Then, exploiting the definition of $T$ given in \eqref{eq:T},  we aim to upper-bound the tail probabilities of $T$ with a combination of probabilities associated to the  events $\{A_i (t)- U_i(t)<0\}$,  $i\in\{1,2\}$, for different  time instant $t$.   However,  in doing so, the following difficulty arises.  $\bold{A}(t)$ depends on $\bold{U}(t)$, which itself depends on the selected strategy and on the past trajectory $\bold{A}(\tau)-\bold{U}(\tau)$ for $\tau<t$.
This because, as already mentioned in Remark \ref{rem1},
{whatever strategy is considered, we} can choose a node in the community $G_i$ at time $\tau$ only if $A_i(\tau-1)- U_i(\tau-1)>0$.
We refer to this constraint as {\em feasibility} constraint.

By (\ref{prop:equiv}), we can choose whatever strategy is convenient (among those 
that are feasible, i.e., satisfy the feasibility constraint),  indeed the choice of a strategy has no impact on the final number of active nodes.  
A first crucial step in our proofs consists in identifying such a strategy.
In the attempt to balance the number of active and unused nodes in the two communities, 
 a possible candidate is the {\em max-strategy}, according to which, at time step $1\leq t\leq T$, 
one chooses the community with the maximum number of active and unused nodes $A_i(t-1)-U_i(t-1)$. The main drawback of this strategy is that the analysis of the corresponding processes $\bold{A}(t)$,  $\bold{U}(t)$, $t\leq T$, appears prohibitive due to its complex correlation structure.
To circumvent this difficulty, we introduce a {\em hybrid} variant of the max-strategy defined above,  according to which,   at time $t$,  the community $G_1$ is selected if and only if
\[
\lim_{n\to \infty} \frac{E[A_1(t)-U_1(t)\mid\bold{U}(t)=\bold{u}(t)]}{g_1}
\ge \lim_{n\to \infty} \frac{E[A_2(t)-U_2(t)\mid\bold{U}(t)=\bold{u}(t)]}{g_2},
\]
i.e., at time $t$ we select the community with the largest
asymptotic normalized expected number of active and unused nodes. 

We go on selecting communities according to this rule up to a random time $T'$, $T'\le T$,  defined as the first time at which the feasibility constraint prevents us from further using our deterministic policy.  For every  time $t\in (T',T]$, instead, we select communities according to an arbitrary {\em feasible} strategy, such as the max-strategy. 
The reason why the {\em hybrid max-strategy} simplifies
the analysis of the bootstrap percolation process is that
up to time $t\leq T'$,  the process $\bold{U}(t)$ is  deterministic, 
with the mapping $t\mapsto(U_1(t)/g_1,U_2(t)/g_2)$ describing a particular 
well determined curve in $\oD$.    As a result, the characterization of $P(A_i(t)-U_i(t)<0)$ becomes extremely simple,
since  it can be reduced to the tail probability of  binomial random variables.  Then we can easily bound from above the probability
$P(A_i(t)-U_i(t)<0)$ by  using the concentration inequalities reported  in Appendix \ref{Penrose}, provided that we  are able to characterize the  average asymptotic dynamics of $E[A_i(t)-U_i(t)]$. We emphasize that, by so doing,  we obtain exponential bounds.
Moreover we wish to point out that  the asymptotic analysis of the average dynamics of the hybrid max-strategy permits us to identify three regimes, which are shown to be equivalent to
$({\bold{Sub}})$,  $({\bold{Sup}})$ and   $({\bold{Crit}})$.

At last we recall that the interested reader can find the extension to the case of SBMs with $k>2$ communities 
in \cite{TGL2}. 
While the stochastic analysis can be carried out following the same lines as for the case $k=2$, the identification of a suitable deterministic strategy is not straightforward. We report in Appendix \ref{subsec:kcomm} a brief discussion of the main issues arising in the case of $k > 2$ communities.}

\section{Proofs}\label{sec:proofs}

\subsection{Preliminaries}\label{subsec:03112021}

We start by introducing the asymptotic normalized mean number of active and unused nodes.
For \textcolor{black}{$t\in[n]\cup\{0\}$}  and $i\in\{1,2\}$, we set
\begin{equation}\label{eq:R}
R_i(\bold{u}(t)):=E[A_i(t)-U_i(t)\,|\,\bold{U}(t)=\bold{u}(t)]=a_i+(n_i-a_i)b_i(\bold{u}(t))-u_i(t).
\end{equation}
Hereon, for $\bold{x}:=(x_1,x_2)\in [0,\infty)^2$, we set
\[
\lfloor \bold{x}g\rfloor:=(\lfloor x_1 g_1\rfloor,\lfloor x_2 g_2\rfloor),
\]
where $\lfloor x\rfloor$ denotes the greatest integer less than or equal to $x\in\mathbb R$. The following lemmas hold.

\begin{lemma}\label{le:bas2}
Assume \eqref{eq:a1a3}, \eqref{eq:gratiog}, \eqref{eq:hyp2bis}, \eqref{eq:trivial} and
let $i\in\{1,2\}$ be fixed. Then
\begin{equation}\label{eq:Btob}
\lim_{n\to\infty}\frac{R_i(\lfloor\bold{x}g\rfloor)}{g_i}=\rho_i(\bold x),
\quad\text{$\forall$ $\bold{x}\in [0,\infty)^2$}
\end{equation}
\end{lemma}

\begin{lemma}\label{le:June15}
Assume \eqref{eq:a1a3}, \eqref{eq:gratiog}, \eqref{eq:hyp2bis}, \eqref{eq:trivial}, and
let $\mathcal W$ be a compact subset of {$(0,\infty)^2$}. Then
\begin{equation*}
\sup_{\bold{x}\in\mathcal{W}}\Big|\frac{R_i(\lfloor \bold{x}g\rfloor)}{g_i}-\rho_i(\bold x)\Big|\to 0,
\quad\text{$\forall$ $i=1,2$.}
\end{equation*}
\end{lemma}

\begin{lemma}\label{le:aspi}
Assume \eqref{eq:a1a3}, \eqref{eq:gratiog} and \eqref{eq:hyp2bis}. Then, for any $\bold{x}\in [0,\infty)^2$, $i\in\{1,2\}$ and $j\in\{1,2\}\setminus\{i\}$,
\begin{align*}
b_i(\lfloor\bold{x}g\rfloor)&=
\left(1+O\left(\bold{1}\{x_i>0\}(\lfloor x_i g_i\rfloor p_i+(\lfloor x_i g_i\rfloor)^{-1})+
\bold{1}\{x_j>0\}(\lfloor x_j g_j\rfloor q+(\lfloor x_j g_j\rfloor)^{-1})\right)\right)\nonumber\\
&\qquad\qquad
\times\left(\lfloor x_i g_i\rfloor p_i+\lfloor x_j g_j\rfloor q\right)^{r}/r!.
\end{align*}
\end{lemma}

We postpone the proofs of Lemmas \ref{le:bas2}, \ref{le:June15} and \ref{le:aspi}, which are technical, but conceptually rather straightforward, to 
Appendix \ref{applemmas}. 


\subsection{Equivalent formulations of $(\bold{Sub})$, $(\bold{Crit})$ and $(\bold{Sup})$}\label{sect:equivk2}

Throughout this subsection we assume $(C)$ and  \eqref{eq:alfaorder} with $\alpha_1\leq 1$.
We consider the curve
\begin{equation}  \label{Drho-def}
\mathcal{D}_{\bm{\rho}}:=\{\bold x\in\mathcal{D}:\,\,\rho_1(\bold x)=\rho_2(\bold x)\}
\end{equation}
and the conditions:\\
\noindent$({\mathcal Sub})$: $\min_{\bold x\in\oD_{\bm\oldb}}\oldb_1(\bold{x})<0$, \qquad  $({\mathcal Crit})$: $\min_{\bold x\in\oD_{\bm\oldb}}\oldb_1(\bold{x})=0$,\qquad $({\mathcal Sup})$: $\min_{\bold x\in\oD_{\bm\oldb}}\oldb_1(\bold x)>0$.\\

\noindent Note that $\mathcal{D}_{\bm{\rho}}$ is graphically represented by the purple curve in  Figures \ref{fig:fig_sub}, \ref{fig:fig_crit} and \ref{fig:fig_sup}. The following proposition holds.

\begin{proposition}\label{le:equivcondboldcall}
Under the assumption $(C)$ with $\alpha_1\leq 1$, we have that the conditions $(\bold{Sub})$, $(\bold{Crit})$ and $(\bold{Sup})$
are equivalent to $({\mathcal Sub})$, $({\mathcal Crit})$ and $({\mathcal Sup})$, respectively.
\end{proposition}

The proof of this proposition exploits the following lemma. 

\begin{lemma}\label{le:Drho}
Assume $(C)$ with $\alpha_1\leq 1$. Then:\\
\noindent$(i)$ $\mathcal{D}_{\bm{\rho}}$ is the graph of a strictly increasing function of class $C^1$, say
$\zeta(\cdot)$, with domain $[x_1^{(0)},x_1^{(1)}]$, where $x_1^{(0)}$ is the unique solution of the equation
\[
\rho_1(x_1,0)-\rho_2(x_1,0)=0,\quad x_1\in \textcolor{black}{ (0,r/(r-1))}
\]
and $x_1^{(1)}$ is the unique point in $(0,r/(r-1))$ such that
\[
\{(x_1^{(1)},\zeta(x_1^{(1)}))\}=\widetilde{\mathcal{D}}\cap\mathcal{D}_{\bm{\rho}},
\]
where
\[
\widetilde{\oD}:=\left\{\bold{x}\in\oD:\,\,\max\{x_1+\chi_{12}x_2, x_2+\chi_{21}x_1\}=\frac{r}{r-1}\right\}.
\]
\noindent$(ii)$ $\rho_1(\bold{x}^{(0)})=\rho_2(\bold{x}^{(0)})>0$, where $\bold{x}^{(0)}:=(x_1^{(0)},0)$.\\
\noindent$(iii)$ $\widetilde{\mathcal{E}}_1$ is the graph of a strictly increasing and strictly concave function of class $C^2$, say
$\zeta_1(\cdot)$, with domain $[y_1^{(0)} ,y_1^{(1)}]$ and $\zeta_1( y_1^{(0)})=0$. Here $ y_1^{(0)}$ is the smallest 
positive solution of the equation $ \alpha_1-x_1+ r^{-1}(1- r^{-1})^{r-1}x_1^r=0$,
and $y_1^{(1)}$ is the unique point on $(0,r/(r-1))$ such that
\[
\{(y_1^{(1)},\zeta_1(y_1^{(1)}))\}=\widetilde{\mathcal{D}}\cap\widetilde{\mathcal{E}}_1.
\]
\noindent$(iv)$ $\widetilde{\mathcal{E}}_2$ is the graph of a strictly increasing and strictly convex function of class $C^2$, say
$\zeta_2(\cdot)$, with domain $[0,y_1^{(2)}]$ and $\zeta_2(0)=y_2^{(0)}$. Here $y_2^{(0)}$  is the smallest positive solution of the equation $ \alpha_2-x_2+ r^{-1}(1- r^{-1})^{r-1}x_2^r=0$,
and $y_1^{(2)}$ is the unique point on $(0,r/(r-1))$ such that
\[
\{(y_1^{(2)},\zeta_2(y_1^{(2)})\}=\widetilde{\mathcal{D}}\cap\widetilde{\mathcal{E}}_2.
\]
\end{lemma}
{
	\textcolor{black}{
 Having established the above lemma, we define
\[
\mathcal{Z}:=\widetilde{\mathcal{E}}_1\cap\widetilde{\mathcal{E}}_2,
\]
i.e., $\mathcal{Z}$ is the set of the zeros of both $\rho_1(\cdot, \cdot)$ and 
$\rho_2(\cdot, \cdot)$, which necessarily lie in $\mathcal D_{\rho}$.}
}
{ \color{black}Under the assumption $(C)$,  by Lemma \ref{le:Drho} (parts $(iii)$ and $(iv)$) 
we have that:
\begin{equation} \label{primaimplicaz}
 \mathcal{Z}=\widetilde{\mathcal{E}}_1 \cap  \widetilde{\mathcal{E}}_2=\emptyset \Leftrightarrow
 \zeta_1(x)<\zeta_2(x),  \; \forall  x\in[y_1^{(0)}, y_1^{(1)} ]\cap[0, y_2^{(1)} ] 
 \Rightarrow  { \mathcal{E}}_1\cap {\mathcal{E}}_2=\emptyset \Leftrightarrow (\bold{Sup}).  
\end{equation}
Hence 
\begin{equation}\label{secondaimplicaz}
(\bold{Sub})  \Leftrightarrow \overset\circ{ \mathcal{E}}_1 \cap   \overset\circ{\mathcal{E}}_2 \neq  \emptyset
\Rightarrow { \mathcal{E}}_1\cap {\mathcal{E}}_2\neq \emptyset  \Rightarrow
 \mathcal{Z}=\widetilde{\mathcal{E}}_1\cap\widetilde{\mathcal{E}}_2\neq \emptyset. 
\end{equation} 
Let
\begin{equation}\label{eq:zerostar}
\bold{z}_*=(z_*,\zeta(z_*))\in\widetilde{\mathcal{E}}_1\cap\widetilde{\mathcal{E}}_2
\end{equation}
denote the zero of $\rho_1(\cdot,\cdot)$ and $\rho_2(\cdot,\cdot)$ in $\mathcal{D}_{\bm{\rho}}$ with the smallest first coordinate (which is obviously
strictly positive), and set
\begin{equation}\label{eq:xstar}
x_*:=z_*+\zeta(z_*)(\nu\mu^r)^{1/(r-1)}>0.
\end{equation}

Here, $\zeta(\cdot)$ is the function whose graph is $\mathcal{D}_{\bm{\rho}}$ (see Lemma \ref{le:Drho}$(i)$).

For later purposes, it is important to note that, as immediate consequence of Lemma \ref{le:Drho} (parts $(iii)$ 
and $(iv)$) we have that, under the assumptions $(C)$ and $(\bold{Sub})$,
there exists a right neighborhood of $z_*>0$, say $I_{z_*}^{+}$, such that $\zeta_2(x_1)>\zeta_1(x_1)$ for any $x_1\in I_{z_*}^{+}$, i.e.,
\begin{equation}\label{eq:intzstar}
\rho_1(x_1,\zeta(x_1))=\rho_2(x_1,\zeta(x_1))<0,\quad\forall\,\,x_1\in I_{z_*}^{+}.
\end{equation}
The proofs of Lemma \ref{le:Drho}  and Proposition \ref{le:equivcondboldcall} are reported in 
 Appendix \ref{proof-deterministiche}.

\subsection{Proof of Theorem \ref{teo-blocksub}}

By Remark \ref{re:alfamag1} we necessarily have $\alpha_1\leq 1$. 
Let $x_1^{(0)}, x_1^{(1)}$ be the extreme points of the domain of $\zeta(\cdot)$ (see Lemma \ref{le:Drho}$(i)$),
consider the segment
\[
\mathcal{S}:=\{(x_1,0):\,\,x_1\in [0,x_1^{(0)}]\}
\]
and denote by $\overline{\zeta}(\cdot)$ the function whose graph is given by
$\mathcal{C}:=\mathcal{S}\cup\mathcal{D}_{\bm{\rho}}$, i.e.,
\[
\overline{\zeta}(x_1):=\bold{1}_{[x_1^{(0)},x_1^{(1)}]}(x_1)\zeta(x_1),\quad x_1\in [0,x_1^{(1)}].
\]
We recall that, in our terminology, a strategy is a rule according to which at every time step \textcolor{black}{$t\in[n]$} 
a community is selected, see Remark \ref{re:strategy}.

We proceed by dividing the proof in four steps. Usually, throughout the proof, for ease of notation, we denote by $c>0$ a generic positive constant, by $c(\varepsilon)$ if it depends on $\varepsilon>0$.

\subsubsection*{\bf Step\,\,1:\,\,\textcolor{black}{Identification\,\,of}\,\,a\,\,suitable\,\,strategy}
{\color{black} 
In this section we are going to formally define the  {\em hybrid} variant 
of the max-strategy, which has been introduced informally in Sect.ion \ref{sec:high-level}. We  start from  its initial deterministic component, in correspondence of which the trajectory of the 
normalized number of used nodes in each community follows the curve $(x_1, \zeta(x_1))$, as it 
can be  observed by combining \eqref{eq:R}, Lemma \ref{le:bas2}, \eqref{Drho-def} and Lemma \ref{le:Drho} (i).
Therefore, our first goal is to define the corresponding un-normalized  \lq trajectory' $(w_1(t), w_2(t))$
of  the actual number of nodes to be used by time $t$ in the community $G_i$.	
With this in mind, \textcolor{black}{we first establish a map between the discrete parameter $t$ and
the 
quantity $x_1$ that parametrizes
the curve $(x_1, \zeta(x_1))$}. In particular, we set}
\begin{equation}\label{eq:tau}
v(x_1):=\lfloor x_1 g_1\rfloor+\lfloor\overline{\zeta}(x_1)g_2\rfloor,\quad x_1\in [0,x_1^{(1)}].
\end{equation}
Note that $v([0,x_1^{(1)}])$ is a subset of 
$[n]\cup\{0\}$, say
$v([0,x_1^{(1)}])=\{t_0,t_1,\ldots,t_{m+1}\}$. Without loss of generality, we assume
$t_0:=v(0)=0<t_1<\ldots<t_m<t_{m+1}:=v(x_1^{(1)})$.
We consider the right-continuous generalized inverse function of $v$:
\begin{equation}\label{eq:taumeno1}
v^{-1}(t_{s}):=\inf\{x_1\in [0,x_1^{(1)}]:\,\,v(x_1)\geq t_{s}\},\quad\text{$s=0,\ldots,m+1$}.
\end{equation} 
Finally,  we set 
\begin{equation}\label{eq:zetasullimag}
w_1(t_{s}):=\lfloor v^{-1}(t_{s})g_1\rfloor,\quad w_2(t_{s}):=\lfloor\overline{\zeta}(v^{-1}(t_{s}))g_2\rfloor,
\quad\text{$s=0,\ldots,m+1$.},
\end{equation}
Now, to conclude our construction,
we extend the definition of $w_i(\cdot)$, $i=1,2$, to the set $(0,v(x_1^{(1)}) )\cap(\mathbb N\cup\{0\})$, by
interpolating their values in $v([0,x_1^{(1)}])$ as follows. We note that by construction
\[
w_i(t_{s+1})-w_i(t_{s})\in\{0,1\},\quad\text{$i=1,2$, $s=0,\ldots,m$}
\]
and
\[
\sum_{i=1}^{2}w_i(t_s)=t_s,\quad\text{$s=0,\ldots,m+1$.}
\]
So, for any $s=0,\ldots,m+1$, $t_{s+1}-t_s\in\{1,2\}$. Consequently, for any $t\in ((0,v(x_1^{(1)}))\cap\mathbb N)\setminus v([0,x_1^{(1)}])$,
there exists $t_s\in v([0,x_1^{(1)}])$, for some $s\in\{0,\ldots,m\}$, such that $t=t_s+1$ and
$t+1=t_{s+1}$. For such a $t$, we define
\begin{equation}\label{eq:zetaext}
w_1(t):=\lfloor v^{-1}(t_{s})g_1\rfloor+1
\quad\text{and}\quad
w_2(t):=\lfloor\overline{\zeta}(v^{-1}(t_{s}))g_2\rfloor.
\end{equation}
Note that, by construction,
\[
\sum_{i=1}^{2}w_i(t)=t,\quad\text{$\forall$ $t\in\{0,\ldots,v(x_1^{(1)})\}$.}
\]
\textcolor{black}{Finally, we need to determine the conditions under which the deterministic strategy we are defining can be successfully employed. To this purpose we define the stopping time}
{\color{black}
\begin{equation}
T':=\min\{1\leq t\leq v(x_1^{(1)}):\,\,A_i(t-1)<w_i(t),
\quad\text{for some $1\leq i\leq 2$}\}.\label{eq:tnprimo}
\end{equation}
}At time step $t\geq 0$, we choose the community $G_i$, $i=1,2$, if and only if $C_i(t)=1$, where
\begin{equation}\label{eq:policydet}
C_i(0):=0,\quad C_i(t):=w_i(t)-w_i(t-1),\quad\text{$1\leq t<\textcolor{black}{ T'}$}
\end{equation}
and
\textcolor{black}{
\begin{equation}\label{eq:strategy-add}
C_1(t)=\bold{1}\{A_1(t-1)-U_1(t-1) \ge A_{2}(t-1)-U_2(t-1)\},\quad C_2(t)=1-C_1(t),\quad \text{$\textcolor{black}{T'}\leq t\leq T$}
\end{equation}
}
In words, the chosen strategy is deterministic and equal to \eqref{eq:policydet} as long as possible.
Indeed, $T'$ is the first time at which the deterministic strategy \eqref{eq:policydet} can not be employed
because of the lack of usable and active nodes. 
Note that setting $\bm{w}(t):=(w_1(t),w_2(t))$, we have
{\color{black}
\begin{equation}\label{eq:policyadd1}
\bold{U}(t)=\bm{w}(t),\,\,\forall\,\, t< T'
\end{equation}
}
 \textcolor{black}{The} above strategy is well-defined, indeed, by construction,
at each time step $t\leq T'$, there exists only one index $i\in\{1,2\}$
such that $C_i(t)=1$ ($C_j(t)=0$ for $j\neq i$), and by \eqref{eq:policyadd1} we have
$T \geq T'$ with $T-1=|\mathcal{G}|$.  As already mentioned, we extend the process for $T\le t\le n$ by adopting an arbitrary \lq\lq unfeasible" strategy.
\textcolor{black}{The choice of the strategy employed for $t\ge T'$ has no impact on $T'$.}
\subsubsection*{\bf Step\,\,2:\,\,outline\,\,of\,\,the\,\,proof.}
It can be easily seen that
\[
\lim_{\delta\to 0}\lim_{n\to\infty}v(z_*\pm\delta)/v(z_*)=1.
\]
Therefore, for any $\varepsilon>0$, there exist $\delta_\varepsilon>0$ and $n_{\varepsilon}\in\mathbb N$ such that 
for any $n\geq n_{\varepsilon}$,
it holds $v(z_*+\delta_\varepsilon)<(1+\varepsilon)v(z_*)$ and
$v(z_*-\delta_\varepsilon)>(1-\varepsilon)v(z_*)$. So, for an arbitrarily fixed $\varepsilon>0$ and any $n\ge n_\varepsilon$ 
{\color{black}
\begin{align}
\{|\mathcal{G}|/v(z_*)-1|>\varepsilon\}&=\{|\mathcal{G}|>(1+\varepsilon)v(z_*))\}\cup\{|\mathcal{G}|<(1-\varepsilon)v(z_*)\}\nonumber\\
&\subseteq\{|\mathcal{G}|>v(z_*+\delta_\varepsilon)\}\cup\{ |\mathcal{G}| < v(z_*-\delta_\varepsilon)
 \subseteq\{|\mathcal{G}|\ge v(z_*+\delta_\varepsilon)\}\cup\{T'\le v(z_*-\delta_\varepsilon)\}.\nonumber
\end{align} 
Since $v(z_*)/g_1\to x_*$, }
\textcolor{black} {the claim then follows if we prove that, for any $\delta>0$ small enough there exists a positive constant $c(\delta)>0$ such that}
{
\begin{align}
&P(\textcolor{black}{T'}\textcolor{black}{\le}v(z_*-\delta))=O(\mathrm{e}^{-c(\delta)g_1})\label{eq:Iconvinprob}\\
&P(|\mathcal{G}|\ge v(z_*+\delta))=O(\mathrm{e}^{-c(\delta)g_1}).  \label{eq:IIconvinprob}
\end{align}
}
\vspace{-0.5 cm}
\subsubsection*{\bf Step\,\,3:\,\,proof\,\,of\,\,\eqref{eq:Iconvinprob}} 
We divide the proof of \eqref{eq:Iconvinprob} in three parts.
In Step 3.1  
we prove the inequality
\begin{equation}\label{eq:T1T2}
P(\textcolor{black}{T'\le}v(z_*-\delta)) \textcolor{black}{<\mathfrak{T}_1+\mathfrak{T}_{2}},
\end{equation}
where
{\color{black}
\begin{align}
\mathfrak{T}_{1}:=\sum_{s=0}^{v(x_1^{(0)})-1}
P(\mathrm{Bin}(n_1-a_1,b_1((w_1(s),0)))<w_1(s+1)-a_1)\label{eq:hug0n}\\
\mathfrak{T}_{2}:=\sum_{s=v(x_1^{(0)})}^{v(z_*-\delta)}
\sum_{i=1}^{2}P(\mathrm{Bin}(n_i-a_i,b_i(\bm{w}(s)))<w_i(s+1)-a_i).
\label{eq:uso5bisn}
\end{align}
}
In Step 3.2 we prove
\textcolor{black}{
\begin{equation}\label{eq:T1zero}
\mathfrak{T}_{1}=O(\mathrm{e}^{-c g_1})
\end{equation}
by applying the concentration inequalities  reported in 
Appendix  \ref{Penrose} to every addend in 
\eqref{eq:hug0n}}. In Step 3.3 we prove
\textcolor{black}{
\begin{equation}\label{eq:T2zero}
\mathfrak{T}_{2}=O(\mathrm{e}^{-c(\delta)g_1})
\end{equation}}
 \textcolor{black}{
again by applying the concentration inequalities reported in 
	Appendix \ref{Penrose} 
	to every addend in \eqref{eq:uso5bisn}.}

\noindent {\it { \bf Step\,\,3.1:}\,\,proof\,\,of\,\,\eqref{eq:T1T2}.}\\
Hereon, we set $\bold{A}(t):=(A_1(t),A_2(t))$ and, for two vectors $\bold{y}=(y_1,y_2)$ and
$\bold{y}'=(y'_1,y'_2)$, we write $\bold{y}\geq\bold{y}'$ if $y_i\geq y'_i$, $i=1,2$.
From \eqref{eq:tnprimo} and \eqref{eq:policyadd1} we get
\textcolor{black}{
\begin{align}
\{T'>t\}&=
\{\bold{A}(s)\geq\bm{w}(s+1)\,\,
\forall\; 0\le s\le t-1\}\subseteq \{\bold{U}(s)= \bm{w}(s)\,\,\forall \; 0\le s\le t\},\label{eq:defTgretat}
\end{align}
which yields}
{\color{black}
\begin{align}
\{T'=t\}&=\{\bold{A}(t-1)<\bm{w}(t), \, \bold{A}(s)\ge\bm{w}(s+1)\,\,
\forall\; 0\le s\le t-2\}\\ & \subseteq \{ \bold{A}(t-1)<\bm{w}(t),  \, \bold{U}(s)=\bm{w}(s)\,\,\forall \;0 \le s\le t-1\}.
\label{eq:defTequal}	
\end{align}
}
{\color{black} 
Therefore 
\begin{align}
P(T'\le t) & =P\left(\bigcup_{1\le s\le t}\{T'=s\}\right)= \sum_{s=1}^ t  P(T'=s)\nonumber \\
&\le \sum_{s=0}^{t-1} P(\bold{A}(s)<\bm{w}(s+1), \bold{U}(h)=\bm{w}(h)\,\,\forall 0\le h\le s ) \\
&\le \sum_{s=0}^{t-1} P(\bold{A}(s)<\bm{w}(s+1), \bold{U}(s)=\bm{w}(s))
\label{added-EL2}
\end{align}
}
Consequently,
{\color{black}
\begin{align}
&P(\textcolor{black}{T'\le} v(z_{*}-\delta))
 \textcolor{black}{\le \sum_{s=0}^{v(z_{*}-\delta)}}P(\bold{A}(s)\textcolor{black}{<}\bm{w}(s+1)\,\mid\,\bold{U}(s)=\bm{w}(s))\nonumber\\
&= \sum_{s=0}^{v(x_1^{(0)})-1}\Biggl(P(S_1(s)+a_1-w_1(s+1)<0\,\mid\,\bold{U}(s)=(w_1(s),0))\Biggr)\nonumber\\
&\qquad\qquad\qquad+
\textcolor{black}{\sum_{s=v(x_1^{(0)})}^{v(z_{*}-\delta)-1}}
\sum_{i=1}^{2}P(S_i(s)+a_i-w_i(s+1)<0\,\mid\,\bold{U}(s)=\bm{w}(s)),\label{eq:disprod}
\end{align}}
where we used the fact that $w_2(s)=0$ for $s=1,\ldots,v(x_1^{(0)})$. The inequality \eqref{eq:T1T2} 
follows from \eqref{eq:disprod}, 
noticing that \eqref{eq:bin} yields
{\color{black}
\begin{align*}
\sum_{s=0}^{v(x_1^{(0)})-1} P(S_1(s)+a_1-w_1(s+1)<0\,\mid\,\bold{U}(s)=(w_1(s),0))=
\mathfrak{T}_{1}
\end{align*}
and
\begin{align*}
\sum_{s=v(x_1^{(0)})}^{v(z_{*}-\delta)-1}
\sum_{i=1}^{2}P(S_i(s)+a_i-w_i(s+1)<0\,\mid\,\bold{U}(s)=\bm{w}(s))=
\mathfrak{T}_2.
\end{align*}
}
\noindent {\it {\bf Step\,\,3.2}:\,\,proof\,\,of\,\,\eqref{eq:T1zero}.} \\
We first note that, since $w_2(s)=0$ and therefore $w_1(s)=s$  for $s < v(x_1^{(0)})$, we have
{\color{black}
\begin{align*}
\mathfrak{T}_{1}
&  = \sum_{s=0}^{\min(a_1, v(x_1^{(0)}))-1} P\Big(\mathrm{Bin}(n_1-a_1,b_1(s,0) )<s+1-a_1 \Big) \\
&+  \sum_{s=  \min(a_1, v(x_1^{(0)}))}^{v(x_1^{(0)})-1} P\Big(\mathrm{Bin}(n_1-a_1,b_1(s,0) )<s+1-a_1 \Big),
\end{align*} 
with the convention that the second addend on the right hand side is null when $\min(a_1, v(x_1^{(0)}))=v(x_1^{(0)})$.
Now, note that by construction    
\[
\sum_{s=1}^{ \min(a_1, v(x_1^{(0)}))-1} P\Big(\mathrm{Bin}(n_1-a_1,b_1(s,0) )<s+1-a_1 \Big)=0.
\] 
Therefore, $\mathfrak{T}_{1}$ is not null only if $a_1 <v(x_1^{(0)})$,  
and so
\[
\mathfrak{T}_{1}\leq\sum_{s= a_1}^{v(x_1^{(0)})-1} P\Big(\mathrm{Bin}(n_1-a_1,b_1(s,0) )<s+1-a_1 \Big).
\]

Now we are going to bound each addend of the sum in the right-hand side  by using the inequality
\textcolor{black}{\eqref{Penrose1.6}} in  Appendix \ref{Penrose}.  
{\color{black}
For any $s\in\{a_1, \ldots,v(x_1^{(0)})-1\}$, we have 
 $v^{-1}(s)=s/g_1= w_1(s)/g_1$.  
 Moreover observe that, since $a_1/g_1\to\alpha_1<x_1^{(0)}$, 
 for $n$ sufficiently large $a_1/g_1\in  [\alpha_1/2, x_1^{(0)} ]$. Similarly, since $v^{-1}(v(x_1^{(0)})-1)= (v(x_1^{(0)})-1) /g_1= 
 (\lfloor x_1^{(0)} g_1 \rfloor -1) /g_1 \uparrow   x_1^{(0)} $, for sufficiently large $n$  we have  $ v^{-1}(v(x_1^{(0)})-1) \in  [\alpha_1/2, x_1^{(0)} ]$.  Then, as an immediate consequence of the monotonicity of the involved functions,  for $n$ sufficiently large,  let us say 
   \mbox{$n>n'$},  we have  that  $v^{-1}(s) \in [\alpha_1/2, x_1^{(0)} ]$  for any  $s\in\{a_1, \ldots,v(x_1^{(0)})-1\}$.
Hence we can apply  Lemma
\ref{le:June15}  and conclude that, for  $n>n''$  ($n''$ not depending on $s$ and not smaller  than $n'$):
\begin{align*}
	R_1((s,0)) & = R_1((w_1(s),0))=R_1((\lfloor v^{-1}(s)g_1\rfloor,0))=R_1(( v^{-1}(s)g_1,0))\\& >\frac{1}{2}(\rho_1(v^{-1}(s),0))g_1
 	\ge  \inf_{x\in[\alpha_1/2, x_1^{(0)} ] } \frac{1}{2}(\rho_1(x,0))g_1
	=\frac{1}{2}\rho_1(\bold{x}^{(0)})g_1,\nonumber
\end{align*}
where for the latter equation we have used the property that   
 $\rho_1(\cdot,0)$ is strictly decreasing on \mbox{$(0,r/(r-1))$}. In conclusion, we have, for all $n>n''$: 
\begin{align}
	E[\mathrm{Bin}(n_1-a_1,b_1(s,0)))]
=R_1((s,0))-a_1+s
\geq \frac{1}{2} \rho_1(\bold{x}^{(0)})g_1+s-a_1
>s+1-a_1.\nonumber
\end{align}
}
Therefore, 
{we can apply   inequality  \textcolor{black}{
\eqref{Penrose1.6}  in Appendix}  \ref{Penrose}.  
Note that, since the mapping $x\mapsto x/(x+y)$, for fixed $y>0$, 
is strictly increasing on $(0,\infty)$,  and the function
$H$ defined by \eqref{eq:H} in Appendix \ref{Penrose} is decreasing on $(0,1)$, for every $n>n''$,} we have
{\color{black}
\begin{align}
&P\Biggl(\mathrm{Bin}(n_1-a_1,b_1((s,0))< s+1-a_1\Biggr)=
P\Biggl(\mathrm{Bin}(n_1-a_1,b_1((s,0))\leq s-a_1\Biggr)\nonumber\\
&\quad\quad
\leq
\exp\Biggl( -\frac{1}{2}\rho_1(\bold{x}^{(0)})g_1
H\left(\frac{s-a_1}{\frac{1}{2}\rho_1(\bold{x}^{(0)})g_1+s-a_1}\right)\Biggr)
\leq\exp\Biggl(-\frac{1}{2}\rho_1(\bold{x}^{(0)})g_1
H\left(\frac{x_1^{(0)}}{\frac{1}{2}\rho_1(\bold{x}^{(0)}) +x_1^{(0)}}\right)\Biggr).\nonumber
\end{align}
}
{\color{black} In conclusion, defined $c':= \frac{1}{2}\rho_1(\bold{x}^{(0)})H\left(\frac{x_1^{(0)}}{\frac{1}{2}\rho_1(\bold{x}^{(0)}) +x_1^{(0)}}\right)$ and $c=\frac{1}{2}c'$,
for every $n\geq n''$ we have
\[
\mathfrak{T}_{1}<  x_1^{(0)} g_1 \mathrm{e}^{-c' g_1}=O( \mathrm{e}^{-c g_1}),
\]
}
which yields \eqref{eq:T1zero}.\\

\noindent {\it {\bf Step\,\,3.3:}\,\,proof\,\,of\,\,\eqref{eq:T2zero}.} \\
\textcolor{black}{ To prove \eqref{eq:T2zero} we follow the same lines as in the proof of \eqref{eq:T1zero}.
As first step we show that  for a sufficiently large $n$ (independently from $s$) $(w_1(s)/g_1,w_2(s)/g_2)$  is contained in a properly compact
set $\mathcal{C}_{z_*-\delta,\varepsilon_0}$ satisfying the following property:
\[
\mathcal{C}_{z_*-\delta,\varepsilon_0}\subset\{\bold{x}\in\mathcal{D}:\,\,\rho_1(\bold{x})>0,\rho_2(\bold{x})>0\}.
\]
Then we bound $\mathfrak{T}_{2}$ as follows:
\begin{align}
\mathfrak{T}_{2}&\leq
\sum_{s=v(x_1^{(0)})}^{v(z_*-\delta)-1}
\sum_{i=1}^{2}P(\mathrm{Bin}(n_i-a_i,b_i(\bm{w}(s))/g_1<(w_i(s)-a_i+1)/g_i )\nonumber\\
&\leq
{v(z_*-\delta)}( \mathfrak{s_1}(\varepsilon)+\mathfrak{s_2}(\varepsilon)),\label{eq:lowerbdsup}
\end{align}
for any $\varepsilon>0$, where 
\begin{equation}\label{eq:mathfraks}
\mathfrak{s_i}(\varepsilon):= \sup_{\bold{x}\in\mathcal{C}_{z_*-\delta,\varepsilon_0}}
P\Biggl(\mathrm{Bin}(n_i-a_i,b_i(\lfloor\bold{x}g\rfloor))/g_i<
x_i-\alpha_i+\varepsilon\Biggr),\quad i=1,2.
\end{equation}
Then, exploiting Lemma \ref{le:June15} and the inequality \textcolor{black}{\eqref{Penrose1.6} in Appendix \ref{Penrose}}, we are going to show that {there exists $\varepsilon=\varepsilon(\delta)$}
such that}
\begin{align}
\mathfrak{s_i}(\varepsilon(\delta))=O(\mathrm{e}^{-c(\delta)g_1}),\quad\text{$i=1,2$.} \label{added-EL}
\end{align}
Then   \eqref{eq:T2zero} immediately follows.}\\
\\
\noindent \textcolor{black}{\it Step 3.3.1:  Definition of $\mathcal{C}_{z_*-\delta,\varepsilon_0}$.} \\
Let
$\mathcal{C}_{x}$ be the graph of the function $\overline{\zeta}(\cdot)$
restricted to $[x_1^{(0)},x]$, for an arbitrary $x\in (x_1^{(0)},x_1^{(1)}]$. 
Clearly, for any $x$, $\mathcal{C}_{x}$ is a compact set of $\mathbb{R}^2$. Using Lemma
\ref{le:Drho}, it is easily seen that, for any $\delta$, 
\[
\mathcal{C}_{z_*-\delta}\subset\mathcal{C}_{z_*}\cap\{\bold{x}\in\mathcal{D}:\,\,\rho_1(\bold{x})>0,\rho_2(\bold{x})>0\}.
\]
For $\varepsilon>0$, let $\mathcal{C}_{z_*-\delta,\varepsilon}$ be the $\varepsilon$-thickening of $\mathcal{C}_{z_*-\delta}$, i.e.,
\begin{equation*}
\mathcal{C}_{z_*-\delta,\varepsilon}:=\{\bold x\in\mathbb R^2:\,\,\mathrm{dist}(\bold x,\mathcal{C}_{z_*-\delta})\leq\varepsilon\},
\end{equation*}
where, for $\mathcal{B}\subset\mathbb R^2$,
\[
\mathrm{dist}(\bold x,\mathcal B):=\inf\{\|\bold x-\bold y\|:\,\,\bold y\in\mathcal B\}
\]
and $\|\cdot\|$ is the Euclidean norm. By the regularity properties of the functions $\rho_i$, $i=1,2$,
easily follows that there exists $\varepsilon_0>0$ small enough such that
\[
\mathcal{C}_{z_*-\delta,\varepsilon_0}\subset\{\bold{x}\in\mathcal{D}:\,\,\rho_1(\bold{x})>0,\rho_2(\bold{x})>0\}.
\]
\\
\noindent \textcolor{black}{\it Step 3.3.2:  proof of the relation  $(w_1(s)/g_1,w_2(s)/g_2) \in \mathcal{C}_{z_*-\delta,\varepsilon_0}$.} \\
We are going to show that there exists a positive integer $n_{\varepsilon_0}$ (not depending on $s$) such that
$
(w_1(s)/g_1,w_2(s)/g_2)\in\mathcal{C}_{z_*-\delta,\varepsilon_0}$ for any $n>n_{\varepsilon_0}$ and
$v(x_1^{(0)})\leq s\leq v(z_*-\delta)-1$. Indeed,  given an arbitrary $n$, for any 
 $v(x_1^{(0)})\leq s\leq v(z_*-\delta)-1$, we have
\[
\frac{\lfloor v^{-1}(s)g_1\rfloor}
{g_1}\leq\frac{w_1(s)}{g_1}\leq\frac{\lfloor v^{-1}(s)g_1\rfloor+1}{g_1}
\]
and
\[
\frac{\lfloor\overline{\zeta}(v^{-1}(s))g_2\rfloor}{g_2}\leq\frac{w_2(s)}{g_2}\leq\frac{\lfloor\overline{\zeta}(v^{-1}(s))g_2\rfloor+1}{g_2}
\]
These relations imply
\[
\Big|\frac{w_1(s)}{g_1}-v^{-1}(s)\Big|\leq 1/g_1
\]
and
\[
\Big|\frac{w_2(s)}{g_2}-\overline{\zeta}(v^{-1}(s))\Big|\leq 1/g_2
\]
for any 
$v(x_1^{(0)})\leq s\leq v(z_{*}-\delta)-1$. Therefore
we can select $n_{\varepsilon_0}$ such that
\[
\|(w_1(s)/g_1,w_2(s)/g_2)-(v^{-1}(s),\overline{\zeta}(v^{-1}(s)))\|\leq\varepsilon_0,
\]
\textcolor{black}{and since $(v^{-1}(s),\overline{\zeta}(v^{-1}(s))\in \mathcal{C}_{z_*-\delta}$ 
for any $v(x_1^{(0)})\leq s\leq v(z_*-\delta)$,  we deduce that 
$(w_1(s)/g_1,w_2(s)/g_2)\in\mathcal{C}_{z_*-\delta,\varepsilon_0}$,
for any $n>n_{\varepsilon_0}$}. \\
\\
%

\noindent \textcolor{black}{\it Step 3.3.3: proof of \eqref{added-EL}.}\\
We shall show \eqref{added-EL} for $i=1$,  indeed the case $i=2$ can be proved similarly.
{\color{black} Setting  $\epsilon(\delta):=\min_{\bold{x}\in\mathcal{C}_{z_*-\delta,\varepsilon_0}}\rho_1(\bold{x})>0$
and $\varepsilon(\delta):=\frac{1}{4}\epsilon(\delta)$, we have
\begin{align}
& \frac{(n_1-a_1)b_1(\lfloor\bold{x}g\rfloor)}{g_1} >
r^{-1}(1-r^{-1})^{r-1}(x_1+\chi_{12}x_2)^r- \frac{1}{4}\epsilon(\delta)  
=x_1-\alpha_1+\oldb_1(\bold x)- \frac{1}{4}\epsilon(\delta)\nonumber\\
&\geq x_1-\alpha_1+\epsilon(\delta)- \frac{1}{4}\epsilon(\delta)
>x_1-\alpha_1+ \frac{3}{4}\epsilon(\delta)\qquad\text{for all $\bold x\in\mathcal{C}_{z_*-\delta,\varepsilon_0}$.}\label{eq:relationunif}
\end{align}
Therefore, by  concentration inequality  \textcolor{black}{\eqref{Penrose1.6} in Appendix \ref{Penrose}} 
for all $n$ large enough, we have
\begin{align}
&\sup_{\bold{x}\in\mathcal{C}_{z_*-\delta,\varepsilon_0}}P\Biggl(\mathrm{Bin}(n_1-a_1,b_1(\lfloor\bold{x}g\rfloor))
\leq \left(x_1-\alpha_1+ \frac{1}{4}\epsilon(\delta)\right)g_1\Biggr)\nonumber\\
&\leq\sup_{\bold{x}\in\mathcal{C}_{z_*-\delta,\varepsilon_0}}
\exp\Biggl(-(x_1-\alpha_1+ \frac{3}{4}\epsilon(\delta))g_1 H\left(\frac{x_1-\alpha_1+ \frac{1}{4}\epsilon(\delta)}{x_1-\alpha_1+ \frac{3}{4}\epsilon(\delta)}\right)\Biggr)\label{eq:useH01}\\
&=O(\mathrm{e}^{-c(\epsilon(\delta))g_1})= O(\mathrm{e}^{-c(\delta)g_1}),\nonumber
\end{align}
}
where in \eqref{eq:useH01} we used \eqref{eq:relationunif} and the fact that $H$ decreases on $(0,1)$.
\subsubsection*{\bf Step\,\,4:\,\,proof\,\,of\,\,\eqref{eq:IIconvinprob}.}
For $\delta>0$,  define  the random time
{
\[
Q(\delta):=\max\{t:\,\,U_1(t)\leq \mathfrak{z}_1,
\,\,U_2(t)\leq \mathfrak{z}_2\}
\]
where $\mathfrak{z}_1:=\lfloor(z_*+\delta)g_1\rfloor$ and $\mathfrak{z}_2:=\lfloor\overline{\zeta}(z_*+\delta)g_2\rfloor$.}
{
Note that \textcolor{black}{by construction  
 \begin{align*}
 \text{either }& U_{1}(Q(\delta))= \mathfrak{z}_1 \qquad  \text{ and } \qquad
 U_{2}(Q(\delta))\le \mathfrak{z}_2,\\
 \text{or }&  U_{1}(Q(\delta))< \mathfrak{z}_1 \qquad  \text{ and }  \qquad
 U_{2}(Q(\delta))=\mathfrak{z}_2 \quad \text{ almost surely.} 
\end{align*}
}
In other words,  defining,  for $\bold v\in\mathbb N^2$, the sets 
\[
\mathcal{F}_{\bold v}:=\mathcal{F}^{(1)}_{\bold v}\cup\mathcal{F}^{(2)}_{\bold v},
\]
\[
\mathcal{F}_{\bold{v}}^{(1)}:=\{(w_1,w_2)\in \mathbb{N}^2:\,\,w_1=v_1,\,\, w_2\le v_2\},\qquad \mathcal{F}_{\bold{v}}^{(2)}:=\{(w_1,w_2)\in \mathbb{N}^2:\,\,w_2=v_2,\,\, w_1\le v_1\},
\]
the random vector $\bold{U}(Q(\delta))$ (whose components are $U_i(Q(\delta))$, $i=1,2$)  almost surely satisfies
\[
\bold{U}(Q(\delta))\in\mathcal{F}_{(\mathfrak{z}_1,\mathfrak{z}_2)}
\quad\text{with}
\quad
\textcolor{black}{|\mathcal{F}_{(\mathfrak{z}_1,\mathfrak{z}_2)}|=\mathfrak{z}_1+\mathfrak{z}_2+1}.
\]
As immediate consequence, we have that almost surely
\[
Q(\delta)= U_1(Q(\delta))+U_2(Q(\delta))\le\mathfrak{z}_1 +\mathfrak{z}_2= v(z_*+\delta).
\]
}
{Therefore}
\begin{align}
\left\{|\mathcal{G}|\ge v(z_*+\delta)\right\}&\subseteq\bigcap_{i=1}^{2}\bigcap_{t\leq v(z_*+\delta)}
\left\{S_i(t)+a_i-U_i(t)\geq 0\right\}\nonumber\\
&\subseteq\bigcap_{i=1}^{2}\{S_i(Q(\delta))+a_i-U_i(Q(\delta))\geq 0\}\nonumber\\
&=\bigcup_{\bold{u}\in\mathcal{F}_{(\mathfrak{z}_1,\mathfrak{z}_2)}}
\bigcap_{i=1}^{2}\{S_i(Q(\delta))+a_i-U_i(Q(\delta))\geq 0,\bold{U}(Q(\delta))=\bold{u}\},\nonumber
\end{align}
{and so}
\begin{align}
&P(|\mathcal{G}|\ge v(z_*+\delta))\leq\sum_{\bold{u}\in\mathcal{F}_{(\mathfrak{z}_1,\mathfrak{z}_2)}}
P\left(\bigcap_{i=1}^{2}\{S_i(Q(\delta))+a_i-U_i(Q(\delta))\geq 0\}\,\Big|\,\bold{U}(Q(\delta))=\bold{u}\right)\nonumber\\
&\leq\textcolor{black}{(\mathfrak{z}_1+\mathfrak{z}_2+1)}\max_{\bold{u}\in\mathcal{F}_{(\mathfrak{z}_1,\mathfrak{z}_2)}}
P\left(\bigcap_{i=1}^{2}\left\{S_i(u_1+u_2)+a_i-u_i\geq 0\right\}
\,\Big|\,\bold{U}(Q(\delta))=\bold{u}\right)\label{eq:ineqpIAN}\\
&\leq\textcolor{black}{(\mathfrak{z}_1+\mathfrak{z}_2+1)}\max_{1\leq j\leq 2}\max_{\bold{u}\in\mathcal{F}_{(\mathfrak{z}_1,\mathfrak{z}_2)}^{(j)}}
P\left(\bigcap_{i=1}^{2}\left\{S_i(u_1+u_2)+a_i-u_i\geq 0\right\}
\,\Big|\,\bold{U}(Q(\delta))=\bold{u}\right).\label{eq:finlim}
\end{align}

Note that, for fixed $j\in\{1,2\}$ and
$\bold{u}\in\mathcal{F}_{(\mathfrak{z}_1,\mathfrak{z}_2)}^{(j)}$,
\begin{align}
&P\left(\bigcap_{i=1}^{2}\left\{S_i(u_1+u_2)+a_i-u_i\geq 0\right\}
\,\Big|\,\bold{U}(Q(\delta))=\bold{u}\right)
\leq P\left(S_j(u_1+u_2)+a_j\geq\mathfrak{z}_j
\,\Big|\,\bold{U}(Q(\delta))=\bold{u}\right)\nonumber\\
&\qquad\qquad
=P(\mathrm{Bin}(n_j-a_j,b_j(\bold u))\geq\mathfrak{z}_j-a_j)
\leq P(\mathrm{Bin}(n_j-a_j,b_j((\mathfrak{z}_1,\mathfrak{z}_2))\geq\mathfrak{z}_j-a_j),\nonumber
\end{align}
where the latter inequality follows from
 the stochastic ordering properties of the binomial distribution with respect to its arguments.
\textcolor{black}{Note, indeed, that $b_j(\bold u )$ (as defined in \eqref{eq:pi}) is increasing with respect to  the components of $\bold u$.}
Combining this inequality with \eqref{eq:finlim} we have
\[
P(|\mathcal{G}|\ge v(z_*+\delta))\leq\textcolor{black}{(\mathfrak{z}_1+\mathfrak{z}_2+1)}\max_{1\leq j\leq 2}P(\mathrm{Bin}(n_j-a_j,b_j(\mathfrak{z}_1,\mathfrak{z}_2)
\geq\mathfrak{z}_j-a_j).
\]
Since
\[
\mathfrak{z}_1+\mathfrak{z}_2+1\sim v(z_*+\delta)
\]
the claim then follows if we prove that, for an arbitrarily fixed $i\in\{1,2\}$, the quantity
\[
P(\mathrm{Bin}(n_i-a_i,b_i(\mathfrak{z}_1,\mathfrak{z}_2))\geq\mathfrak{z}_i-a_i)
\]
goes to zero exponentially fast with respect to $g_1$.
 For this  we  employ again the concentration inequality \textcolor{black}{\eqref{Penrose1.5} in Appendix \ref{Penrose}}.
  \textcolor{black}{Since ideas and  computations are similar to those  in the proof of \eqref{eq:Iconvinprob}, we skip some details.}
By Lemma \ref{le:bas2} we have
\begin{equation}\label{eq:equiv33}
(n_1-a_1)b_1(\mathfrak{z}_1,\mathfrak{z}_2)\sim
((z_*+\delta)-\alpha_1)g_1+\rho_1(z_*+\delta,\overline{\zeta}(z_*+\delta))g_1.
\end{equation}
and
\begin{equation}\label{eq:equiv33bis}
(n_2-a_2)b_2(\mathfrak{z}_1,\mathfrak{z}_2)\sim
(\overline{\zeta}(z_*+\delta)-\alpha_2)g_2+\rho_2(z_*+\delta,\overline{\zeta}(z_*+\delta))g_2.
\end{equation}
Therefore
\begin{align}
\frac{\mathfrak{z}_1-a_1}{(n_1-a_1)b_1(\mathfrak{z}_1,\mathfrak{z}_2)}
\to\frac{(z_*+\delta)-\alpha_1}
{z_*+\delta-\alpha_1+\rho_1(z_*+\delta,\overline{\zeta}(z_*+\delta))}\label{limite}
\end{align}
and
\begin{align}
\frac{\mathfrak{z}_2-a_2}{(n_2-a_2)b_2(\mathfrak{z}_1,\mathfrak{z}_2)}
\to\frac{\overline{\zeta}(z_*+\delta)-\alpha_2}
{\overline{\zeta}(z_*+\delta)-\alpha_2+
\rho_2(z_*+\delta,\overline{\zeta}(z_*+\delta))}.\label{limitebis}
\end{align}
By {Lemma \ref{le:Drho}}, we have that there exists $\delta_0>0$ such that
\[
\max_{1\leq i\leq 2}\rho_i((z_*+\delta),\overline{\zeta}(z_*+\delta))=\epsilon(\delta)<0,\quad\text{for any $0<\delta\leq\delta_0$.}
\]
Therefore, by 
\textcolor{black}{ \eqref{Penrose1.5} in Appendix \ref{Penrose},}
for all $n$ large enough, we have
\[
P\Biggl(\mathrm{Bin}(n_1-a_1,b_1(\mathfrak{z}_1,\mathfrak{z}_2))
\geq\mathfrak{z}_1-a_1\Biggr)
\leq\exp\Biggl(-(n_1-a_1)b_1(\mathfrak{z}_1,\mathfrak{z}_2)
H\left(\frac{\mathfrak{z}_1-a_1}{(n_1-a_1)b_1(\mathfrak{z}_1,\mathfrak{z}_2)}\right)\Biggr),
\]
where   $H(x):=1-x+x\log x,\quad x>0,\quad H(0)=1$.
 \textcolor{black}{The exponential decay of
\[
P(\mathrm{Bin}(n_1-a_1,b_1(\mathfrak{z}_1,\mathfrak{z}_2))\geq\mathfrak{z}_1-a_1)
\]
easily follows combining this latter inequality with \eqref{eq:equiv33} and \eqref{limite}, and using that $H$ increases on $(1,\infty)$.
Reasoning in the same way,  but using \eqref{eq:equiv33bis} and \eqref{limitebis} in place of
\eqref{eq:equiv33} and \eqref{limite}, respectively, one proves the exponential decay of
$P\Biggl(\mathrm{Bin}(n_2-a_2,b_2(\mathfrak{z}_1,\mathfrak{z}_2))\geq\mathfrak{z}_2-a_2\Biggr)$.
 }

\subsection{Proof of Theorem \ref{teosupcrit}}

We give the detailed proof in the case $\alpha_1\leq 1$. The case $\alpha_1>1$ follows along similar computations and it is briefly outlined in   \textcolor{black}{in} Appendix \ref{appendix-Theo3.3}.

We denote by $\overline{\zeta}_{\mathrm{ext}}(\cdot)$ the function whose graph is
\[
\mathcal{C}_{\mathrm{ext}}:=\mathcal{C}\cup\mathcal{R}_{\theta_0},
\]
where $\mathcal{C}$ is defined at the beginning of the proof of Theorem \ref{teo-blocksub}
and $\mathcal{R}_{\theta_0}$, $\theta_0>0$ arbitrarily fixed, is the straight line
\[
\mathcal{R}_{\theta_0}:=\{\bold{x}\in\mathbb{R}^{2}\setminus\mathcal{D}:\,\,\bold{x}=(x_1,\theta_0(x_1-x_1^{(1)})+\overline{\zeta}(x_1^{(1)})),\,\,x_1\geq x_1^{(1)}\},
\]
i.e.,
\[
\overline{\zeta}_{\mathrm{ext}}(x_1):=
{\bold{1}_{[0,x_1^{(1)}]}(x_1)}
\overline{\zeta}(x_1)+\bold{1}_{(x_1^{(1)},\infty)}(x_1)(\theta_0(x_1-x_1^{(1)})+\overline{\zeta}(x_1^{(1)})).
\]
Similarly to the proof of Theorem \ref{teo-blocksub} (see \eqref{eq:tau}), we set
\[
v(x_1):=\lfloor x_1 g_1\rfloor+\lfloor\overline{\zeta}_{\mathrm{ext}}(x_1)g_2\rfloor,\quad x_1\geq 0
\]
{and note that} $v([0,\infty))=\{t_s\}_{s\in\mathbb N\cup\{0\}}$, for some $t_0:=0<t_1<\ldots<t_m<\ldots$.
We define $v^{-1}(t_s)$, $s\in \mathbb{N}\cup\{0\}$,
similarly to \eqref{eq:taumeno1}, with obvious changes (i.e., with $[0,\infty)$ in place of $[0,x_1^{(1)}]$ and with
$\mathbb{N}\cup\{0\}$ in place of $\{0,\ldots,m+1\}$), $w_i(t_s)$, $i=1,2$, $s\in\mathbb{N}\cup\{0\}$, similarly to \eqref{eq:zetasullimag}, and we extend the definition of $w_{i}(\cdot)$ to any $t\in\mathbb{N}\setminus v([0,\infty))$ similarly to \eqref{eq:zetaext}. We define $T'$ as in \eqref{eq:tnprimo}
(with $\{1,\ldots,v(x_1^{(1)})\}$ replaced by $[n]$) and,
similarly to the proof of Theorem \ref{teo-blocksub}, for $t< T'$,  the strategy $\{C_i(t)\}$ defined by \eqref{eq:policydet} is adopted. \textcolor{black}{ For $T'\le t<T$,  we assume that the system switches to the strategy defined by  \eqref{eq:strategy-add}. 
However we wish to emphasize that  the choice of the strategy employed when $t\ge T'$ is completely irrelevant for the proof, as it will become clear in the next subsection.} 

\textcolor{black}{We proceed by giving an outline of the proof and then by dividing the proof itself in five steps.}
Hereon, for ease of notation, we denote by $c>0$ a generic positive constant, by $c(\varepsilon)$ if it depends on $\varepsilon>0$.
\subsubsection*{\bf Outline\,\,of\,\,the\,\,proof}
Let $\varepsilon\in (0,1)$ be small.
Since $|\mathcal{G}|\geq T'$, we have
\[
P(n-|\mathcal{G}|>\varepsilon n)\leq 1-P(T' \textcolor{black}{\ge} \lceil (1-\varepsilon)n\rceil),
\]
\textcolor{black}{where $\lceil x\rceil$ denotes the smallest integer greater than or equal to $x\in\mathbb R$}, and therefore it suffices to show that
\[
P(T'\textcolor{black}{<} \lceil (1-\varepsilon)n\rceil) \geq O(\mathrm{e}^{-c(\varepsilon)g_1}),
\]
for some positive constant $c(\varepsilon)>0$. 
{
{\color{black}
We have
\[
P(T' < \lceil (1-\varepsilon)n\rceil)=
P(T'< kv(x_1^{(1)})  )+ P( kv(x_1^{(1)})   \le T'< \lfloor p_1^{-1}  \rfloor )+
 P(  \lfloor p_1^{-1}  \rfloor \le T'< \lceil (1-\epsilon)n \rceil)
 \]
{for some constant $k\in \mathbb{N}$.}
Therefore the proof is completed if we show that every term on the right-hand side
vanishes exponentially fast for sufficiently large $n$}.  To this aim,  as first step, we give  a preliminary bound on $P(T' \in [t_a,t_b) )$ 
for some $t_a,t_b \in [n]$ with $t_a<t_b$.
 { \color{black} 
\subsubsection*{\bf  Step\,\,1:\ A useful preliminary bound }
 
Note that 
by \eqref{eq:defTequal} we have
\begin{align*}
 \{ T' \in [t_a,t_b)\}&=\cup_{t=t_a}^{t_b-1}  \{\bold{A}(t-1)<\bm{w}(t), \, \bold{A}(s)\ge\bm{w}(s+1)\,\,
\forall\; 0\le s\le t-2\} \\ 
& \subseteq \cup_{t=t_a}^{t_b-1} \{ \bold{A}(t-1)<\bm{w}(t),  \, \bold{U}(s)=\bm{w}(s)\,\,\forall \;0 \le s\le t-1\}.
\end{align*}
Therefore 
\begin{align}
\{ T' \in{[t_a,t_b)}\}\subseteq \{ \bold{U}(t_a-1)=\bm{w}(t_a-1)  \}. \label{eq:incEL}
\end{align}
Since the  paths  $A_i(\cdot)$ and the functions $w_i(\cdot)$
are non-decreasing, we have
\begin{align}
\{ \bold{A}({t_a}-1)\geq\bm{w}(t_b), \bold{U}(t_a-1)=\bm{w}(t_a-1) \} \cap \{ T' \in [t_a,t_b) \}=\emptyset.
\label{eq:excEL}
\end{align}
Combining \eqref{eq:incEL} and \eqref{eq:excEL}, we have
\begin{align*}
&P(\bold{A}(t_a-1)\geq\bm{w}(t_b), \bold{U}(t_a-1)=\bm{w}(t_a-1))+P(  T' \in [t_a,t_b) )\\
&= P(\bold{A}(t_a-1)\geq\bm{w}(b), \bold{U}(t_a-1)=\bm{w}(t_a-1))+P( T' \in [t_a,t_b),  \bold{U}(t_a-1)=\bm{w}(t_a-1) )\\
& \le P( \bold{U}(t_a-1)=\bm{w}(t_a-1) ),
\end{align*}
which yields
\begin{align}
\hspace{-0.5 cm} & P(  T' \in [t_a,t_b) )\le P(  T' \in [t_a,t_b) \mid \bold{U}(t_a-1)=\bm{w}(t_a-1)) \nonumber \\&
\le 1 -  P(\bold{A}(t_a-1)\geq\bm{w}(t_b)\mid  \bold{U}(t_a-1)=\bm{w}(t_a-1))
\le \sum_{i=1}^2P({A_i}(t_a-1)<w_i(t_b)\mid  \bold{U}(t_a-1)=\bm{w}(t_a-1)).
\label{eq:ineq-EL}
\end{align}
}
\vspace{-1 cm}
\subsubsection*{\bf \color{black} Step\,\,2:\,\,Bounding\,\, $P(T'<kv(x_1^{(1)}))$}
{\color{black}
By \eqref{added-EL2}  we have
 \begin{align}
 P(T'<  kv(x_1^{(1)}) ) &\le \sum_{s=1}^{kv(x_1^{(1)})-1 }P(\bold{A}(s)<\bm{w}(s+1)\,\mid\,\bold{U}(s)=\bm{w}(s)) \nonumber\\
&= \mathfrak{T}_{1}+\mathfrak{T}_{2}
\end{align}
with
\begin{align*}
\hspace{-1 cm} \mathfrak{T}_{1}&:=\sum_{s=1}^{v(x_1^{(0)})-1 }P(\bold{A}(s)<\bm{w}(s+1)\,\mid\,\bold{U}(s)=\bm{w}(s)),\quad \mathfrak{T}_{2}&:=\sum_{s=v(x_1^{(0)})}^{ kv(x_1^{(1)}) -1}P(\bold{A}(s)<\bm{w}(s+1)\,\mid\,\bold{U}(s)=\bm{w}(s)).
\end{align*}  

Now, following the same  lines as in the proof of \eqref{eq:T1zero},  we can easily show that
{ \color{black}
\begin{align}
 \mathfrak{T}_{1}=& \sum_{s= \min(a_1, v(x_1^{(0)}))}^{v(x_1^{(0)})}P\left(S_1(s)<s+1-a_1\,\mid\,\bold{U}(s)=(s,0)\right){=}O(\mathrm{e}^{-c g_1}). \label{primopezzo}
\end{align}
with $c:= \frac{1}{4}\rho_1(\bold{x}^{(0)})H\left(\frac{x_1^{(0)}}{\frac{1}{2}\rho_1(\bold{x}^{(0)}) +x_1^{(0)}}\right).$
}\\
%
}
Instead, to prove that {\color{black} $$\mathfrak{T}_{2}{=}O(\mathrm{e}^{-c g_1}),$$} 
we can follow the same approach as in the proof of \eqref{eq:T2zero}.
Hereon, we skip many details and highlight the main differences. Let $\mathcal{C}_x$ be the graph of the function $\overline{\zeta}_{\mathrm{ext}}(\cdot)$
restricted to $(x_1^{(0)},x)$, $x>x_1^{(0)}$, and, for $\varepsilon>0$, let $\mathcal{C}_{x,\varepsilon}$ be the $\varepsilon$-thickening
of $\mathcal{C}_x$. As in the proof of Theorem \ref{teo-blocksub}, one has that
there exists $\varepsilon_0\in (0,1)$ small enough so that
\[
\mathcal{C}_{kv(x^{(1)}_1),\varepsilon_0}\subset\{\bold{x}\in\mathcal{D}:\,\,\rho_1(\bold x)>0,\rho_2(\bold{x})>0\}
\]
and it can be shown that there exists $n_{\varepsilon_0}$ (not depending on $s$) such that $(w_1(s)/g_1,w_2(s)/g_2)\in\mathcal{C}_{\textcolor{black}{kv(x^{(1)}_1)},\varepsilon_0}$ for any $n>n_{\varepsilon_0}$ and any $v(x_1^{(0)})\leq s\leq \textcolor{black}{kv(x^{(1)}_1)-1}$. 
{\color{black}
By the assumption (${\mathcal Sup}$) it follows that $\min_{\bold{x}\in\mathcal{C}_{kv(x^{(1)}_1),\varepsilon_0}}\rho_1(\bold x)=:\epsilon>0$.
Then, proceeding exactly as in the proof of Theorem \ref{teo-blocksub},  we can show that 
for 
 $n$ large enough ${\mathfrak{T}}_{2}\leq \mathfrak{s}^{c(\epsilon) g_1}$, where $\mathfrak{s}$ is defined
as in \eqref{eq:mathfraks}, with $\mathcal{C}_{kv(x^{(1)}_1),\varepsilon_0}$ in place of $\mathcal{C}_{z_*-\delta,\varepsilon_0}$.}
\subsubsection*{\bf \color{black} Step\,\,3:\,\,Bounding\,\, $P(k v(x_1^{(1)})   \le T'< \lfloor p_1^{-1}  \rfloor )$ }
For $n$ so large that \textcolor{black}{so that $\lfloor p_1^{-1}\rfloor-1>kv(x_1^{(1)})$},
define
\[
l:=\min\{\ell\geq k:\,\,p_1 m_\ell\geq 1\},\quad\text{where $m_\ell:=k^{\ell/k}v(x_1^{(1)}) $}
\]
Since $m_l\geq\lfloor p_1^{-1}\rfloor$, we have
\[
[ kv(x_1^{(1)}),\lfloor p_1^{-1}\rfloor]\cap\mathbb N\subseteq\bigcup_{\ell=k }^{l-1}[m_\ell,m_{\ell+1}]\cap\mathbb N.
\]
{\color{black}
Now
\begin{align}
P({k}\, v(x_1^{(1)})   \le T'< \lfloor p_1^{-1}  \rfloor )&\le \sum_{\ell=k}^{l-1}
P(m_{\ell}\le T'< m_{\ell+1}) \nonumber \\
& \le \sum_{i=1}^2 \sum_{\ell=k}^{l-1}P({A_i}(m_{\ell}-1)<{w_i}(m_{\ell+1})\mid  \bold{U}(m_{\ell}-1)=\bm{w}(m_{\ell}-1)) \nonumber
\end{align}
where in the latter inequality we have employed \eqref{eq:ineq-EL}.
Moreover
\begin{align}
P(A_i(m_{\ell}-1)<w_i(m_{\ell+1})\mid  \bold{U}(m_{\ell}-1)=\bm{w}(m_{\ell}-1)) 
=P\left(\mathrm{Bin}(n_i-a_i, b_i(\bm{w}(m_{\ell}-1))<w_i(m_{\ell+1})-a_i)\right). \label{eq:unionbound}
\end{align}}
Therefore,  choosing $k$ large enough and arguing as in the proof of relation $(59)$ in \cite{TGL}, for any $i\in\{1,2\}$, any $\ell\in\{k,\ldots,l-1\}$
and all $n$ large enough, we get
\begin{align}\label{eq:disSPA}
&P(\mathrm{Bin}(n_i-a_i,b_i(\bm{w}(m_{\ell}+1)))<w_i(m_{\ell+1}+1)-a_i)\nonumber\\
&\;\;\; \leq P(\mathrm{Bin}(n_i ,b_i(\bm{w}(m_{\ell}+1)))<w_i(m_{\ell+1}+1))
\leq\mathrm{e}^{-c_1g_1}\mathrm{e}^{-(\ell-\lceil c\rceil)c_2g_1},
\end{align}
for some positive constants $c_1,c_2>0$. 
Finally,  by \eqref{eq:unionbound}, for all $n$ large enough, we have
\begin{align}
P({k}\, v(x_1^{(1)})   \le T'< \lfloor p_1^{-1}  \rfloor ) &\leq \sum_{i=1}^{2}\sum_{\ell=\lceil c\rceil}^{l-1}
P(\mathrm{Bin}(n_i-a_i,b_i(\bm{w}(m_{\ell}-1)))<w_i(m_{\ell+1})-a_i)
\leq c_3\mathrm{e}^{-c_1 g_1},\nonumber
\end{align}
for some positive constant $c_3>0$.}
\subsubsection*{\bf \color{black} Step 4: Bounding $ P(  \lfloor p_1^{-1}  \rfloor \le T'< \lceil (1-\epsilon)n \rceil)$}
Let $c\in (0,1)$ be a small positive constant such that,
for all $n$ large enough $P(\mathrm{Bin}(\lfloor p_1^{-1},p_1)\geq r)\geq 2c$
(see e.g. the proof of Lemma 8.2 Case 3 p. 26 in \cite{JLTV}). For all $n$ large enough we have
{\color{black}
\begin{align}
P(  \lfloor p_1^{-1}  \rfloor \le T'< \lceil (1-\epsilon)n \rceil) 
=P(  \lfloor p_1^{-1}  \rfloor < T'< \lceil cn \rceil) +P(  \lceil cn \rceil  \le T'< \lceil (1-\epsilon)n \rceil).
\label{eq:13J0}
\end{align}
From {\eqref{eq:ineq-EL}}, we have
\begin{align}
P(  \lfloor p_1^{-1}  \rfloor < T'< \lceil cn \rceil) 
&\le 1-P(\bold{A}(\lfloor p_1^{-1}\rfloor-1)\geq\bm{w}(\lceil cn\rceil)
\,|\,\bold{U}(\lfloor p_1^{-1}\rfloor-1)=\bm{w}(\lfloor p_1^{-1}\rfloor-1))\nonumber\\
&
\leq \sum_{i=1}^{2}P(\mathrm{Bin}(n_i-a_i,b_i(\bm{w}(\lfloor p_1^{-1}\rfloor-1)))<w_i(\lceil cn\rceil)-a_i).\label{eq:13J1}
\end{align}
Similarly, we get
\begin{align}
P(  \lceil cn \rceil  \le T'< \lceil (1-\epsilon)n \rceil)
\leq \sum_{i=1}^{2}P(\mathrm{Bin}(n_i-a_i,b_i(\bm{w}(\lceil cn\rceil-1)))<w_i(\lceil(1-\varepsilon)n\rceil)-a_i).\label{eq:13J3}
\end{align}
The following inequalities are proved in the Step 4 of the proof of Proposition 4.1 in \cite{TGL} and
hold for any $i\in\{1,2\}$ and all $n$ large enough:
\begin{align}
&P(\mathrm{Bin}\Big(n_i-a_i,b_i(\bm{w}(\lfloor p_1^{-1}\rfloor-1)) )<w_i(\lceil cn\rceil)-a_i\Big)
\leq c_1\mathrm{e}^{- c_2 n},\quad\text{for some constants $c_1,c_2>0$}\nonumber\\
&P\Big(\mathrm{Bin}(n_i-a_i,b_i(\bm{w}(\lceil c n\rceil-1)))<w_i(\lceil(1-\varepsilon)n\rceil)-a_i\Big)
\leq\mathrm{e}^{-c'(\varepsilon) g_1},\quad\text{for some constant $c'(\varepsilon)>0$.}\nonumber
\end{align}
Therefore 
$P(  \lfloor p_1^{-1}  \rfloor \le T'< \lceil (1-\epsilon)n \rceil)\le \mathrm{e}^{-c'(\varepsilon) g_1}+ c_1\mathrm{e}^{- c_2 n}.$
}


\newpage
\appendix
	\section{Generalization to the SBM with $k$ communities}\label{subsec:kcomm}

Conceptually the generalization to $k>2$ communities can be carried out along similar lines, 
however a few more significant difficulties must be faced, and some workaround is needed.
First, observe that the assumption $(C)$ can be naturally extended to the case of $k>2$ communities, as well as
the definitions of $\oD$ and $\widetilde{\oD}$.
Hereon, we assume that matrix $\bm{\chi}$ is irreducible.  Note, however, that
this assumption does not affect the generality of our results, given that 
more general cases can be traced back to the irreducible case.
The main difficulty in the case of $k > 2$ communities stems from the fact that
$\oD_{\bm \oldb}=\{ {\bold x}\in \oD: \rho_i( {\bold x})=\rho_j( {\bold x}), \forall \; 1\le i,j\le k \}$ is not anymore guaranteed to be the trace of  {a  curve $\overline{\zeta}: \mathbb{R}\to \mathbb{R}^k$ with non decreasing components (as required by every  trajectory followed by a strategy).} 
Therefore we 
have to properly re-define the \lq\lq normalized'' trajectory to be followed by the first deterministic phase of the strategy.
To this end, we start from   the solution of the following Cauchy problem:
\begin{equation}\label{eq:cauchy}
	\bold{x}'(\xi)={\bm\oldb}(\bold{x}(\xi)),\quad\bold{x}(\xi_0):=\bold{x}_0.
\end{equation}
for an appropriate initial condition ${\bold x}_0 \in \oD$,   which can be proved to be component-wise increasing in $\oD$. Then,
 we properly extend it,  so as to create a component-wise increasing curve that connects the origin with a point in $\widetilde{\oD}$.
Given that,  the proofs of main theorems  proceed along the same lines of the case $k=2$.
Details are reported in \cite{TGL2}. 

\section{Proof\,\,of\,\,Proposition\,\,\ref{prop:equiv}}\label{subsec:prop21}

We first prove
	\begin{equation}\label{eq:firstinclusion}
		\cG\subseteq\mathcal{A}(T).
	\end{equation}
	This is equivalent to prove $\bigcup_{h=0}^{H}\cG_h\subseteq\mathcal{A}(T)$, for any $H\in\N\cup\{0\}$.
	We show this claim by induction over $H$. We clearly have $\cG_0=\mathcal{A}(0)\subseteq\mathcal{A}(T)$. Assume
	\[
	\bigcup_{h=0}^H\cG_h\subseteq \cA(T),\quad\text{for some $H\in\N$.}
	\]
	The inclusion \eqref{eq:firstinclusion} follows if we check $\cG_{H+1}\subseteq\mathcal{A}(T)$.
	Take $v\in\cG_{H+1}$ and, reasoning by contradiction, suppose $v\notin\mathcal{A}(T)$.
	By the definition of the bootstrap percolation process $v\in\cG_{H+1}$ has at least $r$ neighbors in
	$\bigcup_{h=0}^H\cG_h$. This set of active nodes is contained in $\mathcal{U}(T)$ due to the inductive hypothesis and relation $\mathcal{A}(T)=\mathcal{U}(T)$.
	Consequently, by \eqref{eq:Mv},  $M_v(T)\geq r$ and so, by \eqref{eq:AMv}, $v\in\mathcal{A}(T)$, which is a contradiction.
	We now prove
	\begin{equation} \nonumber 
		\cG\supseteq\mathcal{A}(T).
	\end{equation}
	For this it suffices to prove that $\mathcal{A}(t)\subseteq\cG$ for any $0\leq t\leq T$. We denote by $\Delta\mathcal{A}(t)$ the set of nodes that become active exactly at time $t$.
	Reasoning by contradiction, assume that there exists at least a $v\in\Delta\mathcal{A}(t)$, $r\leq t\leq T$, such that $v\notin\cG_h$, for any $h\in\N\cup\{0\}$.
	Then there must exist a minimum time $t_0$ with $r\leq t_0\leq T$ such that $\Delta \mathcal{A}(t_0)\not\subseteq\cG$.
	Since $v\in\Delta\cA(t_0)$, it has at least $r$ neighbors in $\cU(t_0)$. By construction we have
	$\cU(t_0)\subseteq\cA(t_0-1)$ and $\cA(t_0-1)\subseteq\cG$. So $v$ has $r$ neighbors in $\cG$. Therefore
	$v\in\cG$, which is a contradiction.\\

\section{Proofs of Lemmas \ref{le:bas2} and \ref{le:June15} and Lemma \ref{le:aspi}} \label{applemmas}

Since the proofs of Lemmas \ref{le:bas2} and \ref{le:June15} exploit Lemma \ref{le:aspi}
	while  the proof of Lemma \ref{le:aspi} is self-contained, we report first the proof of  Lemma \ref{le:aspi} and then  those of Lemmas \ref{le:bas2} and \ref{le:June15}.\\
\subsection{ Proof\,\,of\,\,Lemma\,\,\ref{le:aspi} } \label{proof-lemma5.3}
For $i\in\{1,2\}$, we have
\[
b_i(\lfloor\bold{x}g\rfloor)=P(\mathrm{Bin}(\lfloor x_i g_i\rfloor,p_i)+\mathrm{Bin}(\lfloor x_j g_j\rfloor,q)\geq r),
\quad\text{$j\in\{1,2\}\setminus\{i\}$.}
\]
If $x_1=x_2=0$, then the claim is obvious. If $x_i>0$ and $x_j=0$, then by the third relation in
\eqref{eq:ptcto0bm}, \eqref{eq:gratiog} and formula (8.1) in \cite{JLTV}, for any $m\in\mathbb{N}\cup\{0\}$,
\begin{align}
	P(\mathrm{Bin}(\lfloor x_i g_i\rfloor,p_i)\geq m)=\frac{(\lfloor x_i g_i\rfloor p_i)^m}{m!}(1+O(\lfloor x_i g_i\rfloor p_i+(\lfloor x_i g_i\rfloor)^{-1}))\label{eq:apprbin1}
\end{align}
and, similarly, if $x_j>0$ and $x_i=0$, then
\begin{align}
	P(\mathrm{Bin}(\lfloor x_j g_j\rfloor,q)\geq m)
	=\frac{(\lfloor x_j g_j\rfloor q)^m}{m!}(1+O(\lfloor x_j g_j\rfloor q+(\lfloor x_j g_j\rfloor)^{-1})).\label{eq:apprbin2}
\end{align}
Relations \eqref{eq:apprbin1} and \eqref{eq:apprbin2} clearly give the claim if exactly one component of $\bold x$ is
equal to zero. Now, assume $x_1,x_2>0$. By the independence of the binomial random variables, we have
\begin{align}
	b_i(\lfloor\bold{x}g\rfloor)
	=\sum_{h=0}^{r}P(\mathrm{Bin}(\lfloor x_i g_i\rfloor,p_i)\geq r-h)P(\mathrm{Bin}(\lfloor x_j g_j\rfloor,q)=h)
	+P(\mathrm{Bin}(\lfloor x_j g_j\rfloor,q)\geq r+1).\nonumber
\end{align}
Combining this relation with \eqref{eq:apprbin1} and \eqref{eq:apprbin2}, we have
\begin{align}
	b_i(\lfloor\bold{x}g\rfloor)&=(1+O(\lfloor x_i g_i\rfloor p_i+(\lfloor x_i g_i\rfloor)^{-1}))
	(1+O(\lfloor x_j g_j\rfloor q+(\lfloor x_j g_j\rfloor)^{-1}))\nonumber\\
	&\quad\quad\quad\times
	\sum_{h=0}^{r}\frac{(\lfloor x_i g_i\rfloor p_i)^{r-h}}{(r-h)!}\frac{(\lfloor x_j g_j\rfloor q)^{h}}{h!}\nonumber\\
	&+\frac{(\lfloor x_j g_j\rfloor q)^{r+1}}{(r+1)!}(1+O(\lfloor x_j g_j\rfloor q+(\lfloor x_j g_j\rfloor)^{-1}))
	\nonumber\\
	&=(1+O(\lfloor x_i g_i\rfloor p_i+(\lfloor x_i g_i\rfloor)^{-1}))
	(1+O(\lfloor x_j g_j\rfloor q+(\lfloor x_j g_j\rfloor)^{-1}))\nonumber\\
	&\quad\quad\quad\times
	\frac{(\lfloor x_i g_i\rfloor p_i+\lfloor x_j g_j\rfloor q)^{r}}{r!}\nonumber\\
	&+\frac{(\lfloor x_j g_j\rfloor q)^{r+1}}{(r+1)!}(1+O(\lfloor x_j g_j\rfloor q+(\lfloor x_j g_j\rfloor)^{-1})),\nonumber
\end{align}
from which the claim easily follows.\\
\subsection{ Proof\,\,of\,\,Lemma\,\,\ref{le:bas2}} \label{proof-lemma5.1}
By the definition of $R_i(\cdot)$, we have
	\begin{align}
		\frac{R_i(\lfloor\bold{x}g\rfloor)}{g_i}&=\frac{1}{g_i}(a_i+(n_i-a_i)b_i(\lfloor\bold{x}g\rfloor)
		-\lfloor x_i g_i\rfloor)\nonumber\\
		&=\left(\frac{a_i}{g_i}-\frac{\lfloor x_i g_{i}\rfloor}{g_i}\right)+\frac{n_i-a_i}{g_i}
		b_i(\lfloor\bold{x}g\rfloor).\label{eq:Brel}
	\end{align}
	By \eqref{eq:trivial} it follows
	\[
	\frac{a_i}{g_i}-\frac{\lfloor x_i g_{i}\rfloor}{g_i}\to\alpha_i-x_i.
	\]
	By the second relation in \eqref{eq:ptcto0bm}, \eqref{eq:trivial}, Lemma \ref{le:aspi} and the definition of $g_i$, for $j\neq i$,
	we have
	\begin{align}
		\frac{n_i-a_i}{g_i}
		b_i(\lfloor\bold{x}g\rfloor)&\sim\frac{n_i}{g_i}(x_j g_j q+x_i g_i p_i)^r/r!\nonumber\\
		&=\frac{n_i(g_i p_i)^r}{(r!)g_i}\left(x_j\frac{g_j q}{g_i p_i}+
		x_i\right)^r\nonumber\\
		&=\frac{n_i p_i^{r}g_i^{r-1}}{r!}\left(x_j\frac{g_j q}{g_i p_i}+x_i\right)^r\nonumber\\
		&=\frac{(1-r^{-1})^{r-1}}{r}\left(x_j\frac{g_j q}{g_i p_i}+
		x_i\right)^r\to
		r^{-1}(1-r^{-1})^{r-1}(x_j\chi_{ij}+x_i)^r.\label{eq:limitdefbas}
	\end{align}
	The claim then follows by taking the limit as $n\to\infty$ in \eqref{eq:Brel}.\\
\subsection{ Proof\,\,of\,\,Lemma\,\,\ref{le:June15}} \label{proof-lemma5.2}
For $i\in\{1,2\}$ and $j\in\{1,2\}\setminus\{i\}$, we have
	\begin{align}
		\sup_{\bold{x}\in \mathcal{W}}\Big|\frac{R_i(\lfloor\bold{x}g\rfloor)}{g_i}-\rho_i(\bold x)\Big|
		&\leq\Big|\frac{a_i}{g_i}-\alpha_i\Big|+\frac{1}{g_i}\nonumber\\
		&\qquad\qquad
		+\sup_{\bold{x}\in \mathcal{W}}\Big|\frac{n_i-a_i}{g_i}b_i(\lfloor\bold{x}g\rfloor)
		-r^{-1}(1-r^{-1})^{r-1}(x_j\chi_{ij}+x_i)^r\Big|,\label{eq:supremumin}
	\end{align}
	and so we only need to prove that the supremum in \eqref{eq:supremumin} tends to zero as $n\to\infty$.
	For ease of notation, throughout this proof we denote by $O(x_i,x_j)$ the quantity
	\[
	O(\lfloor x_ig_i\rfloor p_i+\lfloor x_j g_j\rfloor q+\lfloor x_i g_i\rfloor^{-1}+\lfloor x_j g_j\rfloor^{-1}).
	\]
	By Lemma \ref{le:aspi} we have
	\begin{align}
		&\sup_{\bold x\in\mathcal{W}}\Big|\frac{n_i-a_i}{g_i}b_i(\lfloor\bold{x}g\rfloor)
		-r^{-1}(1-r^{-1})^{r-1}(x_j\chi_{ij}+x_i)^r\Big|\nonumber\\
		&\qquad\qquad
		=\sup_{\bold x\in \mathcal{W}}\Big|\frac{n_i-a_i}{g_i}(1+O(x_i,x_j))
		(\lfloor x_j g_j\rfloor q+\lfloor x_i g_i\rfloor p_i)^{r}/r!
		-r^{-1}(1-r^{-1})^{r-1}(x_j\chi_{ij}+x_i)^r\Big|\nonumber\\
		&\qquad\qquad
		\leq\sup_{\bold x\in \mathcal{W}}\Big|\frac{n_i-a_i}{g_i}O(x_i,x_j)(x_j g_j q+x_i g_i p_i)^{r}/r!\Big|\nonumber\\
		&\qquad\qquad\qquad
		+\sup_{\bold x\in \mathcal{W}}\Big|\frac{n_i-a_i}{g_i}(\lfloor x_j g_j\rfloor q+\lfloor x_i g_i\rfloor p_i)^{r}/r!
		-r^{-1}(1-r^{-1})^{r-1}(x_j\chi_{ij}+x_i)^r\Big|.\label{eq:suprema}
	\end{align}
	We start considering the first supremum in the right-hand side of \eqref{eq:suprema}. For some positive constant $c>0$
	(throughout this proof we denote with the same symbol $c$ different constants), we have
	\begin{align}
		&\sup_{\bold x\in \mathcal{W}}\Big|\frac{n_i-a_i}{g_i}O(x_i,x_j)
		(x_j g_j q+x_i g_i p_i)^{r}/r!\Big|\nonumber\\
		&\qquad\qquad
		\leq c(g_j q+g_ip_i+g_j^{-1}+g_i^{-1})
		\frac{n_i-a_i}{g_i}(g_j q+g_i p_i)^{r}\to 0,\nonumber
	\end{align}
	where the limit follows by \eqref{eq:ptcto0bm}, \eqref{eq:gratiog} and \eqref{eq:limitdefbas}.
	As far as the second supremum in the right-hand side of \eqref{eq:suprema} is concerned, by the definition of $g_i$,
	we have
	\[
	g_i=[(1-r^{-1})^{-(r-1)}/(r-1)!]n_i(g_i p_i)^{r}.
	\]
	Therefore, $r!g_i=c^{-1}n_i(g_i p_i)^{r}$, where $c:=r^{-1}(1-r^{-1})^{r-1}$, and so
	\begin{align}
		&\sup_{\bold x\in \mathcal{W}}\Big|\frac{n_i-a_i}{g_i}
		(r!)^{-1}\left(\lfloor x_j g_j\rfloor q+\lfloor x_i g_i\rfloor p_i\right)^{r}-c(x_j\chi_{ij}+x_i)^r\Big|\nonumber\\
		&\qquad\qquad
		\leq c\frac{n_i-a_i}{n_i}
		\sup_{\bold x\in \mathcal{W}}\Big|\left(\frac{\lfloor x_j g_j\rfloor q+\lfloor x_i g_i\rfloor p_i}{g_i p_i}\right)^{r}-
		\left(\frac{x_j g_j q}{g_i p_i}+x_i\right)^{r}\Big|\nonumber\\
		&\qquad\qquad\qquad\qquad
		+c\sup_{\bold x\in \mathcal{W}}\Big|
		\frac{n_i-a_i}{n_i}\left(\frac{x_j g_j q}{g_i p_i}+x_i\right)^{r}
		-(x_j\chi_{ij}+x_i)^r\Big|\nonumber\\
		&\qquad\qquad
		=c\frac{n_i-a_i}{n_i}
		\sup_{\bold x\in \mathcal{W}}\left[\left(\frac{x_j g_j q}{g_i p_i}+x_i\right)^{r}
		-\left(\frac{\lfloor x_j g_j\rfloor q+\lfloor x_i g_i\rfloor p_i}{g_i p_i}\right)^{r}\right]\nonumber\\
		&\qquad\qquad\qquad\qquad
		+c\sup_{\bold x\in \mathcal{W}}\Big|
		\frac{n_i-a_i}{n_i}\left(\frac{x_j g_j q}{g_i p_i}+x_i\right)^{r}
		-(x_j\chi_{ij}+x_i)^r\Big|\nonumber\\
		&\qquad\qquad
		\leq c\frac{n_i-a_i}{n_i}
		\sup_{\bold x\in \mathcal{W}}\left[\left(
		\frac{\lfloor x_j g_j\rfloor q}{g_i p_i}+x_i
		+\frac{q}{g_i p_i}\right)^{r}-
		\left(\frac{\lfloor x_j g_j\rfloor q}{g_ip_i}
		+\frac{\lfloor x_i g_i\rfloor}{g_i}\right)^{r}
		\right]\nonumber\\
		&\qquad\qquad\qquad\qquad
		+c\sup_{\bold x\in \mathcal{W}}\Big|
		\frac{n_i-a_i}{n_i}\left(\frac{x_j g_j q}{g_i p_i}+x_i\right)^{r}
		-(x_j\chi_{ij}+x_i)^r\Big|\nonumber\\
		&\qquad\qquad
		\leq c\frac{n_i-a_i}{n_i}
		\sup_{\bold x\in \mathcal{W}}\left[\left(
		\frac{\lfloor x_j g_j\rfloor q}{g_i p_i}+x_i
		+\frac{q}{g_i p_i}\right)^{r}-
		\left(\frac{\lfloor x_j g_j\rfloor q}{g_ip_i}
		+\frac{\lfloor x_i g_i\rfloor}{g_i}\right)^{r}
		\right]\nonumber\\
		&\qquad\qquad\qquad\qquad
		+c\frac{n_i-a_i}{n_i}\sup_{\bold x\in \mathcal{W}}\Big|
		\left(\frac{x_j g_j q}{g_i p_i}+x_i\right)^{r}
		-(x_j\chi_{ij}+x_i)^r\Big|\nonumber\\
		&\qquad\qquad\qquad\qquad
		+c\sup_{\bold x\in \mathcal{W}}\Big|\frac{n_i-a_i}{n_i}
		(x_j\chi_{ij}+x_i)^r-(x_j\chi_{ij}+x_i)^r\Big|.\label{eq:supremabis}
	\end{align}
	Note that the latter supremum in the right-hand side of \eqref{eq:supremabis} goes to zero as $n\to\infty$. As far as the other two suprema in the right-hand side of \eqref{eq:supremabis},
	note that, for any $\delta>0$ there exists $n_\delta$ such that for any $n\geq n_\delta$
	\[
	\chi_{ij}-\delta<\frac{g_j q}{g_i p_i}<\chi_{ij}+\delta,\quad j\neq i
	\]
	and
	\[
	\frac{q}{g_i p_i}<\delta.
	\]
	Therefore, for any $\delta$ small enough and for all $n$ large enough,
	\begin{align}
		&\sup_{\bold x\in \mathcal{W}}\Big|
		\left(\frac{x_j g_j q}{g_i p_i}+x_i\right)^{r}
		-(x_j\chi_{ij}+x_i)^r\Big|\nonumber\\
		&\leq\sup_{\bold x\in \mathcal{W}}\bold{1}\left\{
		\frac{x_j g_j q}{g_i p_i}\geq x_j\chi_{ij}\right\}
		\left[\left(\frac{x_j g_j q}{g_i p_i}+x_i\right)^{r}
		-(x_j\chi_{ij}+x_i)^r\right]\nonumber\\
		&\qquad\qquad
		+\sup_{\bold x\in \mathcal{W}}\bold{1}\left\{\frac{x_j g_j q}{g_i p_i}<x_j\chi_{ij}\right\}
		\left[(x_j\chi_{ij}+x_i)^r-\left(\frac{x_j g_j q}{g_i p_i}+x_i\right)^{r}
		\right]\nonumber\\
		&\leq\sup_{\bold x\in \mathcal{W}}|
		(x_j\chi_{ij}+x_i+\delta\sup\mathcal{W})^{r}
		-(x_j\chi_{ij}+x_i)^r|+
		\sup_{\bold x\in \mathcal{W}}|(x_j\chi_{ij}+x_i)^r-
		(x_j\chi_{ij}+x_i-\delta\sup\mathcal{W})^{r}|.\nonumber
	\end{align}
	Note that there exists $\varepsilon>0$ such that $\mathcal W\subset [\varepsilon,\infty)^2$, and so
	\begin{align}
		\sup_{\bold x\in \mathcal{W}}|
		(x_j\chi_{ij}+x_i\pm\delta\sup\mathcal{W})^{r}
		-(x_j\chi_{ij}+x_i)^r|&\leq(\sup \mathcal{W})^r(\chi_{ij}+1)^r\Big|
		\left(1\pm\frac{\delta\sup\mathcal{W}}{\varepsilon(\chi_{ij}+1)}\right)^{r}-1\Big|.\nonumber
	\end{align}
	Therefore, by the arbitrariness of $\delta$ we have
	\begin{equation}\label{eq:sup1zeronew}
		\sup_{\bold x\in \mathcal{W}}\Big|
		\left(\frac{x_j g_j q}{g_i p_i}+x_i\right)^{r}
		-(x_j\chi_{ij}+x_i)^r\Big|\to 0.
	\end{equation}
	For the latter supremum, we note that, for all $n$ large enough,
	\begin{align}
		&\sup_{\bold x\in \mathcal{W}}
		\left[\left(\frac{\lfloor x_j g_j\rfloor q}{g_i p_i}+x_i
		+\frac{q}{g_i p_i}\right)^{r}-
		\left(\frac{\lfloor x_j g_j\rfloor q}{g_i p_i}+\frac{\lfloor x_i g_i\rfloor}{g_i}\right)^{r}
		\right]\nonumber\\
		&\qquad\qquad\qquad
		\leq\sup_{\bold x\in \mathcal{W}}
		\left[\left(\frac{\lfloor x_j g_j\rfloor q}{g_i p_i}+x_i+\delta\right)^{r}-
		\left(\frac{\lfloor x_j g_j\rfloor q}{g_i p_i}+
		\frac{\lfloor x_i g_i\rfloor}{g_i}\right)^{r}
		\right]\nonumber\\
		&\qquad\qquad\qquad
		\leq\sup_{\bold x\in \mathcal{W}}
		\Big|\left(\frac{x_j g_j q}{g_i p_i}+x_i+\delta\right)^{r}-
		(x_j\chi_{ij}+x_i)^{r}\Big|\nonumber\\
		&\qquad\qquad\qquad\qquad\qquad\qquad
		+\sup_{\bold x\in \mathcal{W}}
		\Big|\left(\frac{\lfloor x_j g_j\rfloor q}{g_i p_i}+\frac{\lfloor x_i g_i\rfloor}{g_i}\right)^{r}-(x_j\chi_{ij}+x_i)^{r}
		\Big|,\nonumber
	\end{align}
	and one can check that these two latter suprema go to zero as $n\to\infty$ arguing as in the proof of relation \eqref{eq:sup1zeronew}.
%
%
\section{ Proof\,\,of\,\,Proposition\,\,\ref{le:equivcondboldcall} and Lemma \ref{le:Drho}} \label{proof-deterministiche}
 \subsection{ Proof\,\,of\,\,Proposition\,\,\ref{le:equivcondboldcall}.} \label{proof-equivcondboldcall}
We divide the proof in three steps.\\
\noindent{\it Step\,\,1:\,\, 
	{$({\mathcal Sub})$\,\,and $(\bold{Sub})$\,\,are\,\,equivalent.}}\\ 
Clearly $({\mathcal Sub})$ implies $(\bold{Sub})$. Viceversa, 
by \eqref{secondaimplicaz},  $(\bold{Sub})$ implies
$\mathcal{Z}\neq\emptyset$, and so \eqref{eq:intzstar}. This immediately implies $({\mathcal Sub})$.\\
\noindent{\it Step\,\,2:\,\,$({\mathcal Sup})$\,\,and $(\bold{Sup})$\,\,are\,\,equivalent.}\\
{We first show that $({\mathcal Sup})$ implies $(\bold{Sup})$.
Since $\mathcal{Z}\subset\mathcal{D}_{\bm{\rho}}$,
$({\mathcal Sup})$ implies  $\mathcal{Z}=\emptyset$. 
By \eqref{primaimplicaz}
$\mathcal{Z}=\emptyset$  implies $(\bold{Sup})$.}\\
{We now show that $(\bold{Sup})$ implies $({\mathcal Sup})$.} 
If $\mathcal{Z}=\emptyset$, then $\rho_1(\bold x),\rho_2(\bold x)\neq 0$ for any $\bold x\in\mathcal{D}_{\bm{\rho}}$. So by Lemma \ref{le:Drho}$(i)$,
$\rho_1(x_1,\zeta(x_1))=\rho_2(x_1,\zeta(x_1))\neq 0$, for any $x_1\in [x_1^{(0)},x_1^{(1)}]$.
By Lemma \ref{le:Drho}$(ii)$ we have $\rho_2(\bold{x}^{(0)})>0$, then $\rho_1(x_1,\zeta(x_1))=\rho_2(x_1,\zeta(x_1))>0$, for any $x_1\in [x_1^{(0)},x_1^{(1)}]$, i.e.,
$\rho_1(\bold x)=\rho_2(\bold x)>0$ for any $\bold{x}\in\mathcal{D}_{\bm{\rho}}$, which implies $({\mathcal Sup})$.\\
\noindent{\it Step\,\,3:\,\,$({\mathcal Crit})$\,\,and $(\bold{Crit})$\,\,are\,\,equivalent.}\\
This is a consequence of the previous Steps 1 and 2. Indeed, the two sets of conditions:
$\{({\mathcal Sub})$, $({\mathcal Sup})$, $({\mathcal Crit})\}$ and
$\{(\bold{Sub})$, $(\bold{Sup})$, $(\bold{Crit})\}$ are both disjoint and exhaustive.\\

%
%
\subsection{ Proof\,\,of\,\,Lemma\,\,\ref{le:Drho}}
%
\noindent{\it Proof\,\,of\,\,$(i)$.}\\
Set $\sigma(\bold x):=\rho_1(\bold x)-\rho_2(\bold x)$, $\bold{x}\in\mathcal{D}$. A straightforward computation gives
\begin{align}
\partial_{x_1}\sigma(\bold x)&=-1+(1-r^{-1})^{r-1}(x_1+\chi_{12}x_2)^{r-1}-(1-r^{-1})^{r-1}\chi_{21}(x_2+\chi_{21}x_1)^{r-1}\nonumber\\
&\leq-1+(1-r^{-1})^{r-1}(x_1+\chi_{12}x_2)^{r-1}\nonumber
\end{align}
and
\begin{align}
\partial_{x_2}\sigma(\bold x)&=(1-r^{-1})^{r-1}\chi_{12}(x_1+\chi_{12}x_2)^{r-1}+1-(1-r^{-1})^{r-1}(x_2+\chi_{21}x_1)^{r-1}\nonumber\\
&\geq (1-r^{-1})^{r-1}\chi_{12}(x_1+\chi_{12}x_2)^{r-1}.\nonumber
\end{align}
Therefore
\[
\partial_{x_1}\sigma(\bold{x})<0\quad\text{and}\quad\partial_{x_2}\sigma(\bold{x})>0,\quad\text{for any $\bold{x}\in\overset\circ{\oD}$.}
\]
By Dini's implicit function theorem, for any $(u,v)\in\oD_{{\bm\oldb}}$ there exist a neighborhood of $u$, say $I_{u}$, a neighborhood of $v$, say $I_{v}$, and a function $\varphi:I_{u}\to I_{v}$ of class $C^1$ so that, for any $\bold{x}\in I_{u}\times I_{v}\subset\oD$, $\sigma(x_1,\varphi(x_1))=0$ and
\[
\varphi'(x_1)=-\partial_{x_1}\sigma(x_1,\varphi(x_1))/\partial_{x_2}\sigma(x_1,\varphi(x_1))>0.
\]
So $\oD_{{\bm\oldb}}$ is the graph of a strictly increasing function of class $C^1$, say $x_2=\zeta(x_1)$ with $x_1\in [c,d]$ for some
$c,d\in [0,r/(r-1)]$. Clearly, \textcolor{black}{as immediate consequence of the fact that  of $\zeta(\cdot)$ is strictly increasing,}
$d:=x_1^{(1)}$ where $x_1^{(1)}$ is defined in the statement of the lemma.
As far as $c$ is concerned, consider the functions of $x_1\in [0,r/(r-1)]$:
\[
\rho_{1}(x_1,0)=\alpha_1-x_1+r^{-1}(1-r^{-1})^{r-1}x_1^{r}\quad\text{and}\quad
\rho_{2}(x_1,0)=\alpha_2+r^{-1}(1-r^{-1})^{r-1}\chi_{21}^{r}x_1^{r}.
\]
An easy computation shows that $\rho_1(\cdot,0)$ is strictly decreasing on $(0,r/(r-1))$ and $\rho_2(\cdot,0)$ is strictly increasing
on $(0,r/(r-1))$. Moreover, by \eqref{eq:alfaorder} we have $\rho_1(\bold 0)=\alpha_1\geq\alpha_2=\rho_2(\bold 0)$. Since $\alpha_1\leq 1$, we have
\[
\rho_1(r/(r-1),0)=\alpha_1-1\leq 0<\alpha_2+\chi_{21}^{r}/(r-1)=\rho_2(r/(r-1),0).
\]
Therefore there exists a unique $x_1^{(0)}\in (0,r/(r-1))$ such that $\rho_1(\bold{x}^{(0)})=\rho_2(\bold{x}^{(0)})$, where
$\bold{x}^{(0)}=(x_1^{(0)},0)$. Thus $c=x_1^{(0)}$.\\
\noindent{\it Proof\,\,of\,\,$(ii)$.}\\
Since $\rho_2(\bold 0)=\alpha_2\geq 0$ and $\rho_2(\cdot,0)$ is strictly increasing on $(0,r/(r-1))$ we have $\rho_1(\bold{x}^{(0)})=\rho_2(\bold{x}^{(0)})>0$.\\
\noindent{\it Proof\,\,of\,\,$(iii)$\,\,and\,\,$(iv)$.}\\
\textcolor{black}{Since the proofs of parts $(iii)$ and $(iv)$ are similar and both follow the line of the proof of part $(i)$,  we limit ourselves to sketch the proof of part $(iii)$. }
We have
\[
\partial_{x_1}\rho_1(\bold x)=-1+(1-r^{-1})^{r-1}(x_1+\chi_{12}x_2)^{r-1}<0,\quad\text{for any $\bold{x}\in\overset\circ{\oD}$}
\]
and
\[
\partial_{x_2}\rho_1(\bold x)=(1-r^{-1})^{r-1}\chi_{12}(x_1+\chi_{12}x_2)^{r-1}>0,\quad\text{for any $\bold{x}\in\overset\circ{\oD}$.}
\]
By Dini's implicit function theorem, for any $(u,v)\in\widetilde{\mathcal{E}}_1$ there exist a neighborhood of $u$, say $J_{u}$, a neighborhood of $v$, say $J_{v}$, and a function $\psi:J_{u}\to J_{v}$ of class $C^1$ so that, for any $\bold{x}\in J_{u}\times J_{v}\subset\oD$, $\rho_1(x_1,\psi(x_1))=0$ and
\[
\psi'(x_1)=-\partial_{x_1}\rho_1(x_1,\psi(x_1))/\partial_{x_2}\rho_1(x_1,\psi(x_1))>0.
\]
The function $\psi$ is indeed of class $C^2$ and (as can be checked by a simple computation)
$\psi''<0$. So $\widetilde{\mathcal{E}}_1$ is the graph of
a strictly increasing and strictly concave function of class $C^2$, say $x_2=\zeta_1(x_1)$ with domain given in the statement.
%
\section {Outline\,\,of\,\,the\,\,proof\,\, of Theorem \ref{teosupcrit}\,\,for\,\,$\alpha_1>1$ }\label{appendix-Theo3.3}
	Note that, if $\alpha_1>1$, then $\mathcal{D}_{\bm{\rho}}$ may be empty and so Lemma \ref{le:Drho} can not be applied. To circumvent this difficulty,
	we define $\mathcal{C}_{\mathrm{ext}}$ as
	\[
	\mathcal{C}_{\mathrm{ext}}:=\mathcal{S}\cup\mathcal{R}_{\theta_0},
	\]
	where $\mathcal{S}$  is the segment joining $(0,0)$ to $(\frac{r}{r-1},0)$
	and $\mathcal{R}_{\theta_0}$, $\theta_0>0$ arbitrarily fixed, is the straight line
	\[
	\mathcal{R}_{\theta_0}:=\{\bold{x}\in\mathbb{R}^{2}\setminus\mathcal{D}:\,\,\bold{x}=(x_1,\theta_0(x_1- r/(r-1))),\,\,x_1\geq 
	r/(r-1)\},
	\]
	i.e., $\mathcal{C}_{\mathrm{ext}}$ is the graph of the function
	\[
	\overline{\zeta}_{\mathrm{ext}}(x_1):=
	\bold{1}_{(r/(r-1),\infty)}(x_1)(\theta_0(x_1-r/(r-1))), \qquad x_1>0.
	\]
	Since $\alpha_1>1$, we have $\rho_1({\bf x})>0$, $\forall {\bf x}\in [0,\infty)^2$, and so
	$\rho_1(x_1, \overline{\zeta}_{\mathrm{ext}}(x_1)),\rho_2(x_1, \overline{\zeta}_{\mathrm{ext}}(x_1))>0$, for any $x_1\geq 0$.
	Noticed this, the proof can be carried on along similar lines as for the case $\alpha_1\le 1$.
%
\section{Concentration inequalities for binomial random variables} \label{Penrose}
	Throughout this paper we  exploit  extensively some
  classical deviation bounds for the binomial distribution (see e.g. Lemma 1.1 p. 16 in \cite{P}), which we report here
	for the sake of completeness.
	Let the function $H$ be defined  by
	\begin{equation}\label{eq:H}
		H(x):=1-x+x\log x,\quad x>0,\quad H(0)=1,
	\end{equation}
	and set $\mu:=m q$, for $m\in\mathbb{N}$, $q\in(0,1)$. For any
	$0<k<m$, we have:\\
	\noindent if $k\geq\mu$, then
	\begin{equation}\label{Penrose1.5}
		P(\mathrm{Bin}(m,q)\geq k)\leq\exp\left(-\mu H\left(\frac{k}{\mu}\right)\right);
	\end{equation}
	if $k\leq\mu$, then
	\begin{equation}\label{Penrose1.6}
		P(\mathrm{Bin}(m,q)\leq k)\leq\exp\left(-\mu H\left(\frac{k}{\mu}\right)\right).
	\end{equation}
	if $k\geq\mathrm{e}^2\mu$, then
\begin{equation}\label{Penrose1.7}
P(\mathrm{Bin}(n,p)\geq k)\leq\exp\left(-\left(\frac{k}{2}\right)\log\left(\frac{k}{\mu}\right)\right).
\end{equation}
\end{document}